\documentclass[11pt]{article}

\usepackage{amssymb}
\usepackage{amsmath}
\usepackage{amscd}
\usepackage{latexsym}
\usepackage{amsthm}
\usepackage{eucal}
\usepackage{enumerate}


\newtheorem{thmSec}{Theorem}[subsection]
\newtheorem{lemSec}[thmSec]{Lemma}
\newtheorem{remNonumber}[thmSec]{Remark}
\newtheorem{propSec}[thmSec]{Proposition}
\newtheorem{proposition}[thmSec]{Proposition}
\newtheorem{claim}[thmSec]{Claim}
\newtheorem{corSec}[thmSec]{Corollary}
\newtheorem{example}[thmSec]{Example}
\newtheorem*{question}{Question}

\newtheorem{maintheorem}{Theorem}

\newtheorem{Corollary}[maintheorem]{Corollary}

\newtheorem{fact}[maintheorem]{Fact}

\numberwithin{equation}{subsection}

\newcommand{\trans}{{}^t \!}

\newcommand{\set}[2]{\{{#1}:{#2}\}}
\newcommand{\Mmax}{M^{\max}}

\newcommand{\Nmax}{N^{\max}}
\newcommand{\Pmax}{P^{\max}}
\newcommand{\rarrowsim}{\smash{\mathop{\,\rightarrow\,}\limits
  ^{\lower1.5pt\hbox{$\scriptstyle\sim$}}}}
\newcommand{\larrowsim}{\smash{\mathop{\,\leftarrow\,}\limits
  ^{\lower1.5pt\hbox{$\scriptstyle\sim$}}}}

\newcommand{\nbar}[1]{\overline{n}_{#1}}
 
\newcommand{\kakup}[3]{K^+({#1},{#2};{#3})}

\newcommand{\kaku}{K(\zeta, \zeta')}
\newcommand{\Har}[2]{\mathcal{H}^{#1}(\mathbb{R}^{#2})}

\newcommand{\hatpi}{\hat{\varpi}}
\newcommand{\wL}{L^2(\mathbb R_+, r^{m-2}dr)}
\newcommand{\WL}{L^2(\mathbb R^m, \frac{dx}{|x|})}
\newcommand{\SL}{SL(2, \mathbb R)\widetilde{}}
\newcommand{\Mp}{Sp(m, \mathbb R)\widetilde{}}
\newcommand{\R}{R}
\newcommand{\tilI}{\widetilde{I}}

\date{}

\begin{document}

\title{% 
The inversion formula and holomorphic extension of 
the minimal representation of the conformal group
}

\author{
        Toshiyuki \textsc{Kobayashi}\footnote{Partially supported 
 by Grant-in-Aid for Scientific Research 
18340037,
Japan Society for the Promotion of Science.}
        \ and
        Gen \textsc{Mano} 
\\
\\
  RIMS, Kyoto University, \\
  Sakyo-ku, Kyoto, 606-8502, Japan \\[\smallskipamount]
\small
\texttt{toshi@kurims.kyoto-u.ac.jp};
\texttt{gmano@kurims.kyoto-u.ac.jp}
}

\maketitle

\vspace*{-3ex}
\begin{center}
\textit{Dedicated to Roger Howe on the occasion of his 60th birthday}
\end{center}

\medskip

\begin{abstract}
The minimal representation $\pi$ of the indefinite orthogonal group
$O(m+1,2)$ is realized on the Hilbert space of square integrable functions
on $\mathbb R^m$ with respect to the measure 
$|x|^{-1} dx_1\cdots dx_m$. 
This article gives an
 explicit integral formula 
for the holomorphic extension of $\pi$ to a holomorphic semigroup
of $O(m+3, \mathbb C)$ by means of the  Bessel function.
Taking its `boundary value',
we also find the integral kernel of the `inversion operator'
corresponding to the inversion element on the 
Minkowski space $\mathbb R^{m,1}$.

\vspace{1ex}
{\it Mathematics Subject Classifications (2000)} : %22E30,22E45,33C10,43A80.
Primary 22E30; % (1973-now) Analysis on real and complex Lie groups
Secondary 22E45, %  (1973-now) Representations of Lie and linear algebraic 
                 %             groups over real fields: analytic methods
33C10  %(1991-now), % Bessel and Airy functions, cylinder functions, $_0F_1$
35J10, % (1973-now) Schr\"odinger operator
43A80, % (1973-now) Analysis on other specific Lie groups
43A85, %  (1973-now) Analysis on homogeneous spaces
47D05, % (1973-1990), % Semigroups of operators
51B20. % (1980-now) Minkowski geometries

\vspace{1ex}
\end{abstract}

\tableofcontents

\section{Introduction}

\subsection{%
Semigroup generated by a differential operator $D$
}
\label{Differential operator D}
Consider the differential operator 
\begin{align} \label{1.1.1}
  D:={}
& |x| \Bigl(\frac{\Delta}4-1\Bigr)
\nonumber
\\
= {}
& \frac{1}{4} \Bigl(\sum_{j=1}^m x_j^2\Bigr)^{\frac{1}{2}}
              \Bigl(\sum_{j=1}^m \frac{\partial^2}{\partial x_j^2} - 4\Bigr)
\end{align}
on $\mathbb R^m$. 
A distinguishing feature here is
 that $D$ has the following properties (see Remark \ref{Remark}):

 1) $D$ extends to a self-adjoint operator on $\WL$.

 2) $D$ has only discrete spectra $\set{-(j+\frac{m-1}2)}{j=0, 1,\cdots}$ 
   in $L^2(\mathbb R^m, \frac{dx}{|x|})$.

Therefore, one can define a continuous operator 
$$
e^{tD} = \sum_{k=0}^\infty \frac{t^k}{k!} D^k
$$
 for $\operatorname{Re}t\ge 0$
with the operator norm $e^{-\frac{m-1}2 \operatorname{Re}t}$
satisfying the composition law
$$
  e^{t_1 D}\circ e^{t_2 D}=e^{(t_1+t_2)D}.
$$
Thus, $e^{tD}~(\operatorname{Re}t > 0)$ forms a holomorphic one-parameter
semigroup. 
Besides, the operator $e^{tD}$ is 
self-adjoint if $t$ is real, and is unitary if $t$ is purely imaginary.

We ask: 
\begin{question}
 Find an explicit formula of $e^{tD}$.
\end{question}

In this context, our main results are stated as follows:

\begin{maintheorem}[see Theorem \ref{Theorem-B}]
\label{main theorem}
The holomorphic semigroup $e^{tD}~(\operatorname{Re} t>0)$ is given by
\begin{equation*}
\bigl(e^{tD}u \bigr)(x)= \int_{\mathbb R^m} K^+(x,x';t)u(x')\frac{dx'}{|x'|}
\quad\text{for } u \in L^2(\mathbb{R}^m, \frac{dx}{|x|}),
\end{equation*}
where the integral kernel $K^+(x,x';t)$ is defined by
\begin{align*}
 K^+(x,x';t):={}
& \frac{2e^{-2(|x|+|x'|)\coth {\frac{t}{2}}}}
       {\pi^{\frac{m-1}{2}}\sinh^{m-1}\frac{t}{2}}
  \tilde{I}_{\frac{m-3}{2}}
  \Bigl( \frac{\psi(x,x')}{\sinh\frac{t}{2}} \Bigr)
. 
\end{align*}
Here, we set 
$\tilde{I}_\nu(z) := \left(\frac{z}{2}\right)^{-\nu}I_\nu(z)$ 
($I_\nu$ is the $I$-Bessel function (see \eqref{definition of I-Bessel}) and 
$$
\psi(x,x'):=2 \sqrt{2(|x||x'|+\langle x, x' \rangle)}
                =4|x|^\frac{1}2|x'|^\frac{1}2\cos \frac{\theta}2,
$$ 
where $\theta\equiv \theta(x, x')$ is the Euclidean 
angle between $x$ and $x'$ in $\mathbb R^m$.
\end{maintheorem}

Particularly important is the special value at $t=\pi \sqrt{-1}$. We set
\begin{align*}
  K^+(x,x'):= {}
&\lim_{\varepsilon \downarrow 0}K^+(x, x'; \pi \sqrt{-1}+\varepsilon)
\\
= {}
&  \frac{2^\frac{m-1}2 \psi(x,x')^{-\frac{m-3}2}}
           {\sqrt{-1}^{m-1} \pi^\frac{m-1}2 } 
  J_\frac{m-3}2 (\psi(x, x'))
\\
= {}
& \frac{2}{\sqrt{-1}^{m-1} \pi^{\frac{m-1}{2}}}
  \tilde{J}_{\frac{m-3}{2}} (\psi(x,x')).
\end{align*}
\begin{Corollary}[see Theorem \ref{thm:D} and Corollary \ref{cor:E}]
The unitary operator $e^{\pi\sqrt{-1}D}$ on $\WL$ 
is given by the Hankel-type transform:
\begin{equation*}
 T: \WL \to \WL, \quad 
  u \mapsto \int_{\mathbb R^m} K^+(x, x') u(x') \frac{dx'}{|x'|}.
\end{equation*}
This transform has the following properties:
\begin{alignat*}{2}
  &\text{(Inversion formula)} \quad T^{-1}=(-1)^{m+1} T. \\
  &\text{(Plancherel formula)} \quad \| Tu\|_{L^2(\mathbb R^m, \frac{dx}{|x|})}
  =\| u\|_{L^2(\mathbb R^m, \frac{dx}{|x|})}.
\end{alignat*}
\end{Corollary}

Let us explain the backgrounds and motivation of our question from three different 
viewpoints:

(1) Hermite semigroup and its variants 
(see Subsection \ref{Our operator D and Hermite operator D}).

(2) $\mathfrak{sl}_2$-triple of differential operators on
$\mathbb{R}^m$ 
(see Subsection \ref{The action of SL(2,R) times O(m)}).

(3) Minimal representation of the reductive group $O(m+1, 2)$ 
(see Subsection \ref{Minimal representation as hidden symmetry}).

\subsection{Comparison with the Hermite operator $\mathcal D$}
\label{Our operator D and Hermite operator D}
Let us compare our operator $D$ on
$L^2(\mathbb R^m, \frac{dx}{|x|})$ with 
the well-known operator 
$\mathcal{D}$ on $L^2(\mathbb R^m, dx)$ defined by
\begin{equation}\label{1.2.1}
  \mathcal D:=\frac{1}4 (\Delta-|x|^2). 
\end{equation}
We call this operator the {\it Hermite operator} 
following the terminology of R. Howe and E.-C. Tan \cite{xHoTa}.
Analogously to $D$, 
the Hermite operator $\mathcal D$ satisfies the following properties:

1) $\mathcal D$ extends to a self-adjoint operator on $L^2(\mathbb R^m, dx)$.

 2) $\mathcal D$ has only discrete spectra 
$\set{-\frac{1}2(j+\frac{m}2)}{j=0, 1,\cdots}$ in $L^2(\mathbb R^m, dx)$.
 
We recall from \cite[\S 5.3 and \S 6]{xHo}  that
 $\mathcal D$ gives rise to a holomorphic 
semigroup $e^{t\mathcal D}~(\operatorname{Re}t>0)$ ({\it Hermite semigroup}).
Then the following results may be regarded as a prototype of
 Theorem \ref{main theorem} 
and Corollary B. 

\begin{fact}[{see \cite[\S 5]{xHo} \cite[\S4.1]{xTa}}]
The holomorphic semigroup 
$e^{t\mathcal D}~(\operatorname{Re}t>0)$ is given by 
\begin{equation*}
   (e^{t\mathcal D}u)(x)= \int_{\mathbb R^m} \mathcal K(x,x';t)u(x')dx',
\end{equation*} 
where $\mathcal K$ is the Mehler kernel defined by
\begin{equation}\label{eqn:Mehler}
  \mathcal K(x,x'; t):= \frac{1}{(2\pi \sinh \frac{t}2)^\frac{m}2}
  \exp\bigl( -\frac{1}2 
      \,  {\vphantom{ \begin{pmatrix} x \\ x' \end{pmatrix}}}^t\!\!\!
 \begin{pmatrix} x \\ x' \end{pmatrix}
  A(t) \begin{pmatrix}x \\x'\end{pmatrix} \bigr).
\end{equation}
Here, we set
$$
  A(t):= \frac{1}{\sinh \frac{t}2} \begin{pmatrix}
 \cosh \frac{t}2 I_m & 
 - I_m \\ -I_m & \cosh \frac{t}2 I_m \end{pmatrix}
\in GL(2m,\mathbb{R}).
$$
\end{fact}

In light of the limit formula: 
$$
  \lim_{\varepsilon \downarrow 0} \mathcal K(x,x'; \pi \sqrt{-1}+\varepsilon)
=\frac{1}{(2\pi\sqrt{-1})^\frac{m}2} e^{-\sqrt{-1} \langle x, x' \rangle}, 
$$
the special value of the operator $e^{t\mathcal{D}}$ 
at $t=\pi\sqrt{-1}$ reduces to the (ordinary) Fourier transform:
\begin{fact}\label{fact:D}
The unitary operator 
$e^{\pi\sqrt{-1}\mathcal D}$ on $L^2(\mathbb R^m, dx)$ is nothing
other than the 
Fourier transform $\mathcal{F}$:
\begin{equation*}
  (\mathcal Ff)(x)= \frac{1}{(2\pi\sqrt{-1})^\frac{m}2} \int_{\mathbb R^m} 
  e^{-\sqrt{-1}\langle x, x'\rangle} f(x') dx'.
\end{equation*}
\end{fact}

We shall see in Section \ref{Integral formula for the inversion operator.} 
a group theoretic interpretation of the fact that $T$ is
unitary and $T^4 = \operatorname{id}$ as well as the fact that
 $\mathcal F$ is unitary and $\mathcal F^4=\operatorname{id}$.

\subsection{The action of $\SL \times O(m)$}
\label{The action of SL(2,R) times O(m)}

The self-adjoint operator $D$ defined by (\ref{1.1.1}) arises in the context of the 
$\mathfrak{sl}_2$-triple of differential operators on $\mathbb R^m$ as follows. 
We define
\begin{equation}
\label{eqn:sl2operator}
  \Tilde h=2\sum_{j=1}^m x_j \frac{\partial}{\partial x_j} +m-1, \quad 
  \Tilde e=2\sqrt{-1} |x|, \quad \Tilde f= \frac{\sqrt{-1}}2 |x| \Delta.
\end{equation}
These operators $\Tilde h, \Tilde e$ and $\Tilde f$ are skew self-adjoint 
operators on $L^2(\mathbb R^m, \frac{dx}{|x|})$, 
and satisfy the $\mathfrak{sl}_2$-relation: 
$$
  [\Tilde h, \Tilde e]=2 \Tilde e, \quad 
  [\Tilde h, \Tilde f]=-2 \Tilde f, \quad 
  [\Tilde e, \Tilde f]= \Tilde h.
$$
The operator $D$ has the following expression
\begin{equation}\label{D}
  D=\frac{1}{2\sqrt{-1}}(-\Tilde e+\Tilde f), 
\end{equation}
which means that $\sqrt{-1}D$ corresponds to a generator of
$\mathfrak{so}(2)$ in $\mathfrak{sl}(2,\mathbb{R})$. 
The $\mathfrak{sl}_2$-module $C_0^\infty(\mathbb R^{m}\setminus
\{0\})$ exponentiates to a unitary representation of 
$SL(2, \mathbb R)$ on $L^2(\mathbb R^{m}, \frac{dx}{|x|})$
(see Subsection \ref{K_0 type decomposition} and Lemma 
\ref{lem:SL2G}).

On the other hand, there is a natural unitary representation of the orthogonal group
$O(m)$ on the same space $L^2(\mathbb R^m, \frac{dx}{|x|})$, and the actions 
of $SL(2, \mathbb R)$ and $O(m)$ mutually commute.
Then we have the following discrete and multiplicity-free decomposition
into irreducible representations of $SL(2, \mathbb R) \times O(m)$
(see \cite[Theorem A]{xkors2}): 
$$
  L^2(\mathbb R^m, \frac{dx}{|x|}) \simeq 
  \sideset{}{^\oplus}{\sum}_{j=0}^{\infty}
  \pi_{2j+m-1}^{SL(2, \mathbb R)} \otimes \Har{j}{m}.
$$  
Here, $\Har{j}{m}$ denotes the 
space of harmonic polynomials on $\mathbb R^m$ of degree $j$, and 
$\pi_{b}^{SL(2, \mathbb R)}$ stands for the irreducible unitary lowest weight 
representation of $SL(2, \mathbb R)$ with minimal $K$-type
$\mathbb C_b$ for $b \in \mathbb N_+ =\{1, 2,  \cdots \}$. 
It is the limit of discrete series if $b=1$, and holomorphic
discrete series if $b \ge 2$.

In contrast, the Hermite operator $\mathcal D$ (see (\ref{1.2.1})) arises from the following  
$\mathfrak{sl}_2$-triple:
\begin{equation}
\label{eqn:sl2Weil}
\Tilde h':=\sum_{j=1}^m x_j \frac{\partial}{\partial x_j}+\frac{m}2, \quad 
\Tilde e':=\frac{\sqrt{-1}}2|x|^2, \quad \Tilde f':=\frac{\sqrt{-1}}2 \Delta, 
\end{equation}
where the Hermite operator
 $\mathcal D$ is given by
\begin{equation}
 \mathcal D=\frac{1}{2\sqrt{-1}}(-\Tilde e'+\Tilde f').
\end{equation}
This $\mathfrak{sl}_2$-triple also
gives rise to the commutative actions of the double covering group $\SL$ 
of $SL(2, \mathbb R)$ and 
$O(m)$ on $L^2(\mathbb R^m)$, 
whose irreducible decomposition amounts to (see \cite[Chapter III, Theorem 2.4.4]
{xHoTa})
$$
  L^2(\mathbb R^m, dx) \simeq 
  \sideset{}{^\oplus}{\sum}_{j=0}^{\infty}
  \pi_{j+\frac{m}2}^{\SL} \otimes \Har{j}{m}.
$$  
Here, $\pi_{b}^{\SL}$ stands for the irreducible unitary lowest weight 
representation of $\SL$ with minimal $K$-type
$\mathbb C_b$ for 
$b \in \frac{1}2\mathbb N_+ =\{ \frac{1}2, 1, \frac{3}2, \cdots \}$. 
It is the Weil representation if $b=\frac{1}{2}$, 
and is obtained by the representation of $SL(2, \mathbb R)$ if
$b \in \mathbb{N}_+$.

\subsection{Minimal representation as hidden symmetry}
\label{Minimal representation as hidden symmetry}

The representation of $SL(2, \mathbb R) \times O(m)$
on $\WL$ in Subsection \ref{The action of SL(2,R) times O(m)} extends to the 
irreducible unitary 
representation $\pi_+$ of the double covering group 
$G:=SO_0(m+1, 2)\widetilde{}$ of 
the indefinite orthogonal group 
(see Subsections \ref{K_0 type decomposition} and \ref{Schrodinger}). 
If $m$ is odd,
this representation is well-defined also as a representation
 of $SO_0(m+1, 2)$. 

Similarly, the representation of $\SL \times O(m)$ on $L^2(\mathbb R^m, dx)$
extends to the unitary representation $\varpi$ of the metaplectic
group $G'=\Mp$. 

These groups $G$ and $G'$ may be interpreted as {\it hidden symmetry}
of $SL(2, \mathbb R)\widetilde{} \times O(m)$.
Conversely, the group $\SL \times O(m)$ forms a `dual pair' in each of
the groups $G$ and $G'$.

The unitary representations $\pi_+$ and $\varpi$ are typical examples 
of \lq minimal representations' of reductive Lie groups in the sense 
that the Gelfand--Kirillov dimension attains its minimum among 
infinite dimensional unitary representations or in the sense that its
annihilator is the Joseph ideal in the enveloping algebra.

The unitary 
representation $\pi_+$ may be interpreted
 as the mass-zero spin-zero wave equation,
or as the bound states of the Hydrogen atom (in $m$ space dimensions),
while the representation $\varpi$ is sometimes referred to 
as the oscillator representation or as the (Segal--Shale--)Weil representation.

We shall review the $L^2$-realization
 of the minimal representation $\pi_+$ of 
$SO_0(m+1, 2)\widetilde{}$,
that is, the analog of the Schr\"{o}dinger model
on $L^2(\mathbb R^m, \frac{dx}{|x|})$ in 
Subsection \ref{Schrodinger}. 
See also \cite{xFo, xHo} for a nice introduction to the original 
Schr\"odinger model of the Weil representation of $\Mp$
on $L^2(\mathbb R^m)$.

To be more precise, we take 
$$
  e=\begin{pmatrix} 0  & 1 \\ 0  &  0 \end{pmatrix},\quad 
  f=\begin{pmatrix} 0  & 0 \\ 1  &  0 \end{pmatrix}, \quad 
  h:=[e,f]=\begin{pmatrix} 1 & 0 \\ 0 & -1 \end{pmatrix}
$$
to be a basis of $\mathfrak{sl}(2, \mathbb R)$ and 
define injective Lie algebra homomorphisms
\begin{align*}
&\phi: \mathfrak{sl}(2, \mathbb R)\to \mathfrak{so}(m+1, 2),
\\
&\varphi:\mathfrak{sl}(2, \mathbb R) \to \mathfrak{sp}(m, \mathbb R)
\end{align*}
(see Subsection \ref{The vector Z}) such that
the differential operators \eqref{eqn:sl2operator} and \eqref{eqn:sl2Weil} 
are obtained via $\phi$ and $\varphi$, respectively, that is, 
\begin{equation*}
  d\pi_+(\phi(X))=\Tilde X, \quad 
  d\varpi (\varphi(X))=\Tilde X'
\end{equation*}
holds for $X=e, f, h$. 
Next we set
$$
  z:= \frac{1}{2\sqrt{-1}} (-e+f) =\frac{1}{2\sqrt{-1}}
  \begin{pmatrix}  0 &-1 \\ 1 & 0 \end{pmatrix} \in 
  \sqrt{-1}\mathfrak{sl}(2, \mathbb R).
$$
Then $e^{\sqrt{-1}\mathbb R \phi(z)}$ is the center of the maximal
compact subgroup of $G$ for $m>1$, while $e^{\sqrt{-1}\mathbb R\varphi(z)}$ is 
that of $G'$ (we use the same notations $\phi$ and $\varphi$ for their complex
linear extensions). 

In this context, 
we shall see in Lemma \ref{differential operator d pi} 
and Remark \ref{Remark2}
that the 
differential operators $D$ and $\mathcal D$ are given by
\begin{equation}\label{eqn:DD}
  D=d\pi_+ (\phi(z)),
 \qquad  \mathcal D= d\varpi(\varphi(z)).
\end{equation}

Thanks to these formulas, Theorem A and Corollary B
are also useful in the analysis on minimal representations of 
$SO_0(m+1, 2)\widetilde{}$ as well as $\Mp$ in the following contexts: 

1) {\bf The Gelfand--Gindikin program --- Theorem A}.

The {\it Gelfand--Gindikin program} asks for extending 
a given unitary representation of 
a real semisimple Lie group $G$ to a holomorphic object of 
some complex submanifold 
in its complexification $G_\mathbb C$.
Stanton and Olshanski\v{\i} \cite{xSt, xOl} independently 
gave a general framework of
 the Gelfand--Gindikin program for holomorphic discrete series.
Their abstract results are enriched, for example for $G'=\Mp$, by the explicit formula 
of the Hermite semigroup $e^{t \mathcal D}=\varpi(e^{t\varphi(z)})$ for the 
Weil representation $\varpi$ on $L^2(\mathbb R^m)$ by Howe \cite{xHo}.
Likewise, 
Theorem A gives an explicit formula of the semigroup $e^{tD}=
\pi_+(e^{t\phi(z)})$ for the minimal representation of 
$G=SO_0(m+1, 2)\widetilde{}$. 
Since $e^{t\phi(z)}\in G_{\mathbb{C}}\setminus G$ for $\operatorname{Re}t>0$,
Theorem~A can be interpreted as a descendent of 
the Gelfand--Gindikin program.

2) {\bf The unitary inversion operator --- Corollary B}.

In the Schr\"odinger model of the minimal representation, $G$ acts only on the 
function space $\WL$ but does not act on the underlying geometry 
$\mathbb R^m$ itself.
One may observe this fact
 by the aforementioned formula $\tilde{f}=d\pi_+(\phi(f))$, which 
does not act on functions on $\mathbb{R}^m$ as a vector field but acts as
 a differential operator of second order on $\mathbb R^m$
(see \eqref{eqn:sl2operator}).
To see how $G$ acts on $\WL$, we use the facts that 

 1) $G$ is generated by $\overline{\Pmax}$ and $w_0$.

 2) The $\overline{\Pmax}$ action on $\WL$ is easily described 
(see Subsection \ref{Schrodinger model of the minimal representation.}).

Here, $\overline{\Pmax}$ is a maximal parabolic subgroup of $G$ 
(see Subsection \ref{Notations and definitions.})
and $w_0:=e^{\pi \sqrt{-1}Z}\in G$ sends $\overline{\Pmax}$ 
to its opposite parabolic subgroup $\Pmax$.
Geometrically, $\overline{\Pmax}$ is essentially 
the conformal affine transformation group on the flat standard 
Lorentz manifold $\mathbb R^{m,1}$ (the \textit{Minkowski space}), 
and $w_0 := e^{\pi\sqrt{-1}Z} \in G$ 
acts on $\mathbb R^{m,1}$
as the \lq inversion'
element (see Subsection \ref{sec:6.1}). 

Thus the representation $\pi_+$ of 
$G$ on $L^2(\mathbb R^m, \frac{dx}{|x|})$ would be understood if we 
find an explicit formula for $\pi_+(w_0)$. 
But since the formula $w_0=e^{\pi \sqrt{-1} \phi(z)}$ implies $\pi_+(w_0)=
e^{\pi\sqrt{-1}D}$, Corollary B answers this question. 
This parallels the fact that the Weil representation is generated by 
the (natural) action of the Siegel parabolic subgroup $P_\text{Siegel}$ 
and the Fourier transform $\mathcal F = e^{\pi\sqrt{-1}\mathcal{D}}$
(see Fact~D).

Briefly, we pin down the analogy in the table below.
Howe \cite{xHo} established the left-hand
 side of the table for the oscillator representation $\varpi$
of $\Mp$,
while Theorem A and Corollary B supply the right-hand side 
of the table for the minimal representation $\pi_+$
of $SO_0(m+1, 2)\widetilde{}$.
$$
\begin{array}{c|cc}
     &  
     \mathfrak{sp}(m, \mathbb R)  
     & \mathfrak{so}(m+1,2)  
\\
\hline
\\
    \text{minimal representation}
     &   
     (d\varpi, L^2(\mathbb R^m))
     & (d\pi_+, L^2(\mathbb R^m, \frac{dx}{|x|})) \\
        e   & \frac{\sqrt{-1}}2 |x|^2
     & 2\sqrt{-1}|x|  \\
        f   & \frac{\sqrt{-1}}2 \Delta
     & \frac{\sqrt{-1}}2|x|\Delta  \\
        h   & E_x+\frac{m}2  &
     2E_x+m-1 \\
z=\frac{-e+f}{2\sqrt{-1}} 
     & \mathcal D:= \frac{1}4 (\Delta-|x|^2) & D:= \frac{1}4|x|\Delta-|x| 
      \\
\text{holomorphic semigroup}~ 
e^{tz}  &  \mathcal{K}(x, x'; t) & K(x, x'; t) \\
\text{inversion}~ 
  e^{\pi \sqrt{-1} z} & \text{Fourier transform} &  \text{Hankel transform}  \\
 \text{maximal parabolic subgroup} &  P_\text{Siegel}  & \overline{\Pmax} 
\end{array}
$$

Analogous results to Corollary B were previously known for some
 singular unitary highest weight representations. 
For example, see a paper \cite{xGrKu} by 
Ding, Gross, Kunze and Richards for those of $U(n,n)$. 
Since $SU(2,2)$ is a double covering group of $SO_0(4,2)$, Corollary  
B in the case $m=3$ essentially corresponds to 
\cite[Corollary 7.5]{xGrKu} in the case 
$(n,k)=(2,1)$ in their notation. 
However, our proof based on an analytic continuation (see
Theorem A) is different from theirs.
Also in \cite{xkmano4}, we shall find the inversion operator
$\pi(w_0)$ for the minimal representation of $O(p,q)$
for $p+q$ even, 
and in particular
give yet another proof of Corollary B
for odd $m$. 

We also present explicit integral 
formulas of $\pi(e^{tZ})$ and $\pi(w_0)$ when restricted to radial functions
and alike (see Theorem \ref{Theorem-C}). This yields a group theoretic 
interpretation of some classic formulas of special functions 
of one variable 
including Weber's second exponential integral formula on Bessel function
and the reciprocal and the Parseval-Plancherel 
formula for the Hankel transform.

\vspace{2ex}
This article is organized as follows. After summarizing
the preliminary results on the $L^2$-model of 
the minimal representation $\pi$,
we find explicitly which function arises for describing $K$-types in $\WL$,
and define a holomorphic extension $\pi_+(e^{tZ})$
in Section \ref{Semigroup pi plus e tZ.}. The integral formula
of the \lq radial' part of $\pi_+(e^{tZ})$ is given in Theorem
\ref{Theorem-C}. 
Theorem~A is proved in Section \ref{Main Theorem}
by using the result of Section \ref{Radial part of the semigroup}.
Taking the special value at $t=\pi\sqrt{-1}$,
we obtain the integral formula of $\pi(w_0)$ 
corresponding to the inversion element $w_0$. 
This corresponds to Corollary~B, and is proved in
Section \ref{Integral formula for the inversion operator.}.
Our integral formula for $\pi(w_0)$ enables us to write explicitly 
the action of the whole group $G = SO_0(m+1,2)\widetilde{}$ on
$L^2(\mathbb{R}^m, \frac{dx}{|x|})$.
This is given in
Section \ref{sec:formulaG}. 
For the convenience of the reader,
we collect basic formulas of special functions in a way that we use in
this article.

The main results of the paper were announced in \cite{xkmano1}
with a sketch of the proof.

Notation:
 $\mathbb{N} = \{ 0, 1, 2, \dots \}$, 
$\mathbb{N}_+ = \{1,2,3,\ldots\}$, 
$\mathbb R_{\le 0}=\set{x \in \mathbb R}{x\le 0}$,
$\mathbb R_+ =\set{x \in \mathbb R}{x >0}$.

\section{Preliminary results on the minimal representation of $O(m+1,2)$}
\label{preliminaries}

This section gives a brief review on the known results of the $L^2$-model of the minimal 
representation of $O(m+1,2)$ in a way that we shall use later. 
We shall give an explicit action of the 
maximal parabolic subgroup $\overline{\Pmax}$ and the Lie algebra 
$\mathfrak g$. Furthermore, we state an explicit $K$-type decomposition 
of $L^2(C_+)$ even though the action of $K$ itself is not given
explicitly here (see Section \ref{sec:formulaG} for this).

\subsection{Maximal parabolic subgroup of the conformal group}
\label{Notations and definitions.}

Let $O(m+1,2)$ be the indefinite orthogonal group which preserves 
 the quadratic form 
$x_0^2+\cdots +x_m^2-x_{m+1}^2-x_{m+2}^2$
of signature $(m+1, 2)$. 
We denote by $e_0, \cdots,e_{m+2}$ the standard 
basis of $\mathbb{R}^{m+3}$, and by $E_{ij}$ $(0\le i, j \le m+2)$ the matrix
 unit.
We set

\begin{alignat*}{1}
 \varepsilon_j&:= \begin{cases}
                 1 \quad  &(1\le j \le m),  \\
                 -1 \quad &(j=m+1).
                \end{cases}  
\end{alignat*}
We take the following elements of the Lie algebra 
$\mathfrak o (m+1,2)$: 
\begin{alignat}{1}
 \overline{N}_j &:= E_{j, 0} + E_{j, m+2} - \varepsilon_j E_{0, j} 
            + \varepsilon_{j} E_{m+2, j} 
\quad (1 \le j \le m+1),  \notag \\ \label{definition of E}
 N_j &:= E_{j, 0} - E_{j, m+2} - \varepsilon_j E_{0, j} 
           - \varepsilon_{j} E_{m+2, j}
\quad (1 \le j \le m+1), \\
 E &:= E_{0, m+2} + E_{m+2, 0},   \notag 
\end{alignat}
and define subalgebras of  $\mathfrak o (m+1, 2)$ by 
\begin{alignat*}{1}
 \overline{\mathfrak{n}^\text{max}}:=
   \sum_{j=1}^{m+1}\mathbb{R}\overline{N}_j, \quad
 \mathfrak{n}^\text{max}:=
   \sum_{j=1}^{m+1} \mathbb{R}N_j, \quad
 \mathfrak{a}:=\mathbb{R}E.
\end{alignat*}
Then we define the following subgroups of $O(m+1,2)$:
\begin{alignat*}{2}
\Mmax_+ &:= \set{g \in O(m+1,2)}{g \cdot e_0 = e_0, \ g \cdot
e_{m+2} = e_{m+2}} 
\\
&\ \simeq \   O(m,1),
\\ 
\Mmax &:= \Mmax_+ \cup \{-I_{m+3}\}\cdot 
\Mmax_+ 
\\
&\ \simeq \   O(m,1) \times \mathbb Z_2,
\\
\overline{\Nmax} &:= \exp (\overline{\mathfrak{n}^\text{max}}), \\
\Nmax& := \exp (\mathfrak{n}^\text{max}),  \\
A &:= \exp (\mathfrak{a}).
\end{alignat*}
For $b=(b_1,\dots,b_{m+1})\in\mathbb{R}^{m+1}$,
we set
\begin{align}\label{definition of n}
\nbar{b}
:={}& \exp(\sum_{j=1}^{m+1} b_j \overline{N}_j)
\\
={} & I_{m+3}+\sum_{j=1}^{m+1} b_j \overline{N}_j + \frac{Q(b)}{2}
      (-E_{0,0} -E_{0,m+1} +E_{m+1,0} +E_{m+1,m+1}),
\nonumber
\end{align}
where $Q(b)$ is the quadratic form of signature $(m,1)$ given by
$$
Q(b) := b_1^2 +\cdots+ b_m^2 - b_{m+1}^2.
$$
The Lie group $\overline{\Nmax}$ is abelian, 
and we have an isomorphism of Lie groups:
$$
\mathbb{R}^{m+1} \rarrowsim \overline{\Nmax}, \quad
b \mapsto \nbar{b}.
$$
It is readily seen from \eqref{definition of n} that
\begin{align}
&\nbar{b}(e_0-e_{m+2}) = e_0 - e_{m+2},
\label{eqn:nbn}
\\
&\nbar{b}(e_0+e_{m+2}) = \trans (1-Q(b), 2b_1,\dots,2b_{m+1}, 1+Q(b)).
\label{eqn:nbp}
\end{align}
We also note
\begin{equation}
\label{eqn:a+m}
e^{tE} (e_0 + e_{m+2})
= e^t (e_0 + e_{m+2}).
\end{equation}
The subgroup $\Mmax_+ \overline{\Nmax}$ is isomorphic to the 
semidirect product group $O(m, 1)\ltimes \mathbb{R}^{m+1}$ via the bijection
$ \mathbb{R}^{m+1} \stackrel{\sim}{\to} \overline{\Nmax}$, 
$ b \mapsto \nbar{b}$.
In this context, 
$\Mmax_+ \overline{\Nmax}$ is regarded as the group of
isometries of the Minkowski space $\mathbb{R}^{m, 1}$,
while $O(m+1 ,2)$ is the group of M\"{o}bius transformations on 
$\mathbb{R}^{m, 1}$ preserving its conformal structure.

Next, we define a maximal parabolic subgroup 
$$
 \overline{\Pmax}:= \Mmax A  \overline{\Nmax}.
$$

In our analysis of the minimal representation of $O(m+1,2)$, $\overline{\Pmax}$
plays an analogous role to the Siegel parabolic subgroup of the metaplectic
group $\Mp$  for the Weil representation.

\subsection{$L^2$-model of the minimal representation}
\label{Schrodinger model of the minimal representation.}

We shall briefly review the $L^2$-model of the minimal representation 
of $O(m+1,2)$. Let $C_\pm$ be the forward and 
the backward light cone respectively: 
$$
  C_\pm:= \set{(\zeta_1, \cdots, \zeta_{m+1}) \in \mathbb R^{m+1}}
  {\pm \zeta_{m+1}>0, \quad\zeta_1^2 + \cdots +\zeta_m^2 =\zeta_{m+1}^2},
$$ 
and $C$ be its disjoint union $C_+ \cup C_+$, that is, $C$ is the conical subvariety
with respect to the quadratic form of signature $(m, 1)$:
\begin{equation}
  C= \set{\zeta \in \mathbb R^{m+1} \setminus \{0\}}
                {Q(\zeta)=0}.
\end{equation} 
Note that $\Mmax_+$ acts on $C$ transitively.

The measure $d\mu$ on $C$ is naturally defined to be 
$\delta(Q)$, the generalized function
associated to the quadratic form $Q$ (see \cite[Chapter~III, \S 2]{xGeSh}). 
Then we form a unitary representation $\pi$ of $\overline{\Pmax}$ on 
the Hilbert space $L^2(C, d\mu)\equiv L^2(C)$ as follows: for $\psi \in 
L^2(C)$, 
\begin{alignat}{2}
\label{rep_A}
   (\pi(e^{t E}) \psi)(\zeta)
   &:= e^{-\frac{m-1}2 t} \psi(e^{-t} \zeta)
   && \qquad (t \in \mathbb R), 
\\
\label{rep_M}
   (\pi(m) \psi)(\zeta)
   &:= \psi(\trans{m} \zeta)
   && \qquad (m \in \Mmax_+), 
\\
\label{rep_m_0}
 (\pi(-I_{m+3}) \psi)(\zeta)
   &:= (-1)^{\frac{m-1}2}\psi(\zeta), 
   &&
\\
\label{rep_N}
   (\pi(\nbar{b}) \psi)(\zeta)
   &:= e^{2\sqrt{-1} (b_1 \zeta_1 + \dots + b_{m+1} \zeta_{m+1})} \psi(\zeta)
   && \qquad (b \in \mathbb R^{m+1}).
\end{alignat}
Then $\pi$ is irreducible and unitary as a $\overline \Pmax$-module, and it is 
proved in \cite[Theorem 4.9]{xkors3} that 
the $\overline{\Pmax}$-module $(\pi, L^2(C))$ extends to 
an irreducible unitary representation of $O(m+1,2)$ if $m$ is odd.
We shall denote this representation of $O(m+1,2)$ by the same letter $\pi$.
The direct sum decomposition 
\begin{equation}
  L^2(C)= L^2(C_+) \oplus L^2(C_-)
\end{equation}
yields a branching law $\pi=\pi_+\oplus \pi_-$ with respect to the restriction 
$O(m+1,2) \downarrow SO_0(m+1,2)$, where $SO_0(m+1,2)$ 
is the identity component of $O(m+1,2)$.
The irreducible representations $\pi_+$ and $\pi_-$ of $SO_0(m+1,2)$ 
are contragredient to each other, one is a highest weight module, and the other is a   
lowest weight module.

\subsection{$K$-type decomposition}
\label{K_0 type decomposition}
Let $SO(2)\widetilde{}$ be the double covering of $SO(2)$.
We write $\eta$ for the unique element of $SO(2)\widetilde{}$ of
order two. 
Then, we have an exact sequence:
$$
1 \to \{1,\eta\} \to SO(2)\widetilde{} \to SO(2) \to 1.
$$
Let 
\begin{equation}\label{eqn:Gdouble}
G=SO_0(m+1,2)\widetilde{}
\end{equation}
be the double covering group of $SO_0(m+1, 2)$ 
 characterized as follows:
a maximal compact subgroup $K$ is of the form $K_1 \times K_2
\simeq SO(m+1) \times SO(2)\widetilde{}$ 
 and the kernel of the 
covering map $G \to SO_0(m+1,2)$ is given by
$\{ (1,1), (1, \eta)\}$. 
Likewise, the double covering group $O(m+1,2)\widetilde{}$ of
$O(m+1,2)$ is defined.

If $m$ is odd, the irreducible representation $(\pi_\pm, L^2(C_\pm))$ defined 
in Subsection \ref{Schrodinger model of the minimal representation.} extends to that 
of 
$SO_0(m+1,2)$, and therefore, also that of $G = SO_0(m+1,2)\widetilde{}$ 
as we proved it more
generally for $O(p,q)$ ($p+q:$ even) in \cite{xkors3}.
We shall use the same letter $\pi_+$ to denote the extension to 
$SO_0(m+1,2)$ or $G$. 
If $m$ is even, by \cite{xSa}, the irreducible unitary representation  
$(\pi_\pm, L^2(C_\pm))$ is still well-defined as 
a representation of $G$, whose Lie algebra 
$\mathfrak g=\mathfrak{so}(m+1,2)$ of $G$ 
acts in the same manner as in the case of odd $m$
(see Subsection \ref{Infinitesimal action.}). 

The $K$-type formula of $(\pi_\pm, L^2(C_\pm))$ is given
as follows:
\begin{equation}
\label{K_0 type decomposition 1}
L^2(C_\pm)_K \simeq \bigoplus_{a=0}^\infty
  \Har{a}{m+1} \boxtimes  \mathbb C 
e^{\pm (a+\frac{m-1}2) \sqrt{-1} \theta}.
\end{equation}
(For example, this formula can be read from \cite[\S1.3]{xSa} by
 substituting 
 $d=m-1, p=1$ and $q=0$.) 
Here, $\Har{a}{m+1}$ stands for the representation of $SO(m+1)$ on
 the space of the spherical harmonics 
which is irreducible if $m>1$
(see Subsection \ref{Spherical harmonics.}). 

Likewise, the representation 
$(\pi, L^2(C))$ of $O(m+1,2)\widetilde{}$
decomposes when restricted to its maximal compact subgroup as follows
 (see \cite[Theorem 3.6.1]{xkors1}):
\begin{equation}
\label{K-type decomposition 1}
  L^2(C)_{K} \simeq \bigoplus_{a=0}^\infty
  \Har{a}{m+1} \boxtimes \Har{a+\frac{m-1}2}{2}.
\end{equation}        
$\Har{a+\frac{m-1}2}{2}$ decomposes into $\mathbb C 
e^{(a+\frac{m-1}2) \sqrt{-1} \theta} \oplus \mathbb C 
e^{-(a+\frac{m-1}2) \sqrt{-1} \theta}$ as $SO(2)$-modules
(see Subsection \ref{Spherical harmonics.}).
This corresponds to the decomposition 
$$
  L^2(C)=L^2(C_+)\oplus L^2(C_-),
$$
for which the $K$-type formula is given by (\ref{K_0 type decomposition 1}).

\subsection{Infinitesimal action of the minimal representation}
\label{Infinitesimal action.}

For $\overline{N}_j, N_j (1\le j \le m+1)$ and $E$,
we define linear transformations on the space 
$\mathcal S'(\mathbb R^{m+1})$ of tempered distributions by
\begin{alignat}{2}
\label{d varpi N j 1}
d \hatpi (\overline{N}_j) 
&:= 2 \sqrt{-1} \zeta_j,
\\
\label{d varpi N j 2}
d \hatpi ({N}_j) 
&:= \sqrt{-1} 
\left( -\frac{m+3}{2}  \varepsilon_j \frac{\partial}{\partial \zeta_j} 
- E_\zeta \varepsilon_j \frac{\partial}{\partial \zeta_j} 
                     + \frac12 \zeta_j \square_\zeta\right),
\\
\label{d varphi E}
d \hatpi  (E) &:= -\frac{m+3}{2}-E_\zeta,
\end{alignat}
where we set
\begin{equation*}
  \square_\zeta :=
    \frac{\partial^2}{\partial \zeta_1^2}+\frac{\partial^2}{\partial \zeta_2^2}
   + \dots +\frac{\partial^2}{\partial \zeta_m^2}
   -\frac{\partial^2}{\partial \zeta_{m+1}^2}, \quad  
   E_\zeta:=\sum_{j=1}^{m+1} \zeta_j 
  \frac{\partial}{\partial \zeta_j}.
\end{equation*}

Then, we recall from \cite{xkors3} that this generates 
the infinitesimal action $d\pi$ of the Lie algebra $\mathfrak{so}(m+1,2)$, and we have
the following commutative diagram for any $X \in \mathfrak{so}(m+1,2)$:
\begin{equation}
\label{commutative diagram}
  \begin{matrix}
    ~~~~~~~~L^2(C)_K  &  \stackrel{\imath}{\longrightarrow} 
                                             &  \mathcal S' (\mathbb R^{m+1}) 
    \\
    d\pi(X) \Bigm\downarrow  &                & ~~~~~~\Bigm\downarrow d\hatpi (X)
    \\
    ~~~~~~~~L^2(C)_K  & \stackrel{\imath}{\longrightarrow}  
                                            & \mathcal S' (\mathbb R^{m+1}).
 \end{matrix}
\end{equation}  
Here, $\imath: L^2(C) \to \mathcal S' (\mathbb R^{m+1})$ is given by 
$u(\zeta) \mapsto u(\zeta) \delta(Q)$. This is 
well-defined and injective if $m>1$ (see \cite[\S 3.4]{xkors3}).

\section{Branching law of $\pi_+$}
\label{Semigroup pi plus e tZ.}

The main goal of this section is Proposition \ref{branching-law},
which explicitly describes special functions that arise as $K$-types 
in the `Schr\"odinger model' on $L^2(\mathbb{R}^m,|x|^{-1}dx)$
of the minimal representation of the double covering group $G$ of
$SO_0(m+1,2)$. 

\subsection{Schr\"odinger model of the minimal representation}
\label{Schrodinger}

We have used the variables $\zeta=(\zeta_1, \cdots, \zeta_{m+1})$ for 
the positive cone 
$C_+ \subset \mathbb R^{m+1}$ in Section \ref{preliminaries},
and will use the letter $x=(x_1, \cdots, x_m)$ for the coordinate of $\mathbb R^m$. 
The projection
\begin{equation}\label{jmath}
p : \mathbb R^{m+1} \to \mathbb R^m, \quad
(\zeta_1, \cdots, \zeta_m, \zeta_{m+1}) \mapsto (\zeta_1, \cdots, \zeta_m).
\end{equation}
induces a diffeomorphism from 
$C_+$ onto $\mathbb R^m \setminus \{ 0\}$,
and the measure $d\mu$ on $C_+$ 
is given by $\delta(Q)$, and therefore
is pushed forward to 
$\frac{1}{2|x|} dx =
 \frac{1}{2|x|} dx_1 \cdots dx_m$. 
Thus, we have a unitary isomorphism: 
\begin{equation}
\sqrt{2} p^\ast :
L^2(\mathbb R^m, \frac{dx}{|x|}) \stackrel{\sim}{\to} L^2(C_+).
\end{equation} 
Through this isomorphism,
we can realize the minimal representation of $G$ on 
$L^2(\mathbb{R}^m, \frac{dx}{|x|})$ as well.
Named after the Schr\"{o}dinger model of the Weil representation,
we say this model is the {\it Schr\"odinger model} of the minimal 
representation of $G$. 
We shall work with this model from now on.

\subsection{$K$-finite functions on the forward light cone $C_+$}
 
This section refines the $K$-type decomposition (\ref{K_0 type decomposition 1})
by providing an explicit irreducible decomposition: 
\begin{equation}\label{R_0 decomposition}
  \WL_{K}=\bigoplus_{a=0}^\infty W_a
 =\bigoplus_{a=0}^\infty\bigoplus_{l=0}^a W_{a,l}
\end{equation}
according to the following chain of subgroups: 
$$
\begin{matrix}
  G & \supset   &
  K  &  \supset & \R:=K \cap (\Mmax_+)_0    \\
  \shortparallel  &  &  \shortparallel  &   &  \shortparallel
  &  &    \\
  SO_0(m+1,2)\widetilde{}  & \supset & SO(m+1)\times SO(2)\widetilde{}
  & \supset & SO(m).
\end{matrix} 
$$

Here, the $K$-irreducible subspace $W_a$ of $\WL$ 
and the $\R$-irreducible supspace $W_{a,l}$ of $W_a$ is characterized by
\begin{align}
&W_a \simeq \Har{a}{m+1} \boxtimes \mathbb C e^{(a+\frac{m-1}2)\theta}
\quad
\text{(see (\ref{K_0 type decomposition 1})),}
\label{eqn:Wa}
\\
&W_{a,l} \simeq \Har{l}{m}, 
\label{eqn:Wal}
\end{align} 
as $K$-modules and as $R$-modules, respectively.
Here, we note that the $SO(m)$-module
$\mathcal{H}^l(\mathbb{R}^m)$ occurs exactly once in the
$O(m+1)$-module $\mathcal{H}^a(\mathbb{R}^{m+1})$ 
if $0 \le l \le a$ 
(see Subsection \ref{Spherical harmonics.} (4)).

Proposition \ref{branching-law} will describe the subspace
 $W_{a,l}$ of $\WL$ by means of Laguerre polynomials $L_n^\alpha(x)$ 
(see \eqref{eqn:defLag} for the definition). 
For this,
we set 
\begin{equation}\label{definition of f_a,l}
f_{a,l}(r):= L_{a-l}^{m-2+2l}(4r)r^l e^{-2r}
\ \ (0 \le l \le a),
\end{equation}
and define injective linear maps by
$$
  j_{a,l}: \Har{l}{m} \to C^\infty (\mathbb R^m), 
 \quad \phi(\omega) \mapsto (j_{a,l}\phi)(r\omega):=f_{a,l}(r) \phi(\omega).
$$
Here, we have identified $\mathbb R^m$ with $\mathbb R_+\times S^{m-1}$
by the polar coordinate
\begin{equation}
\label{polar coordinate}
\mathbb R \times S^{m-1} \to \mathbb R^m, \quad (r, \omega) \mapsto r\omega.
\end{equation}

Then, we have: 

\begin{propSec}
\label{branching-law}
{\rm 1)} 
\begin{equation}
\label{branching-law-2}
j_{a,l}\bigl(\Har{l}{m} \bigr) \subset L^1(\mathbb R^m, \frac{dx}{|x|})
\cap \WL.
\end{equation}

{\rm 2)} Furthermore, the image of $j_{a,l}$ coincides with 
the $\R$-type $W_{a,l}$:
\begin{equation*}
  j_{a,l}\bigl( \Har{l}{m} \bigr) =W_{a,l}.
\end{equation*} 
\end{propSec}

\begin{remNonumber}
\label{rem:falcob}
Proposition \ref{branching-law} (2) asserts in particular that 
$j_{a,l}(\Har{l}m)$ and $j_{a',l}(\Har{l}m)$ are orthogonal 
to each other if $a \ne a'$ $(a, a' \ge l)$. 
More than this, it gives a 
representation theoretic proof of the fact that 
$\set{f_{a,l}(r)}{a=l, l+1,\cdots}$ forms a complete orthogonal basis of $\wL$
(see Lemma \ref{Lemma3.1.1} for the normalization).
\end{remNonumber}

\begin{remNonumber}
The indefinite orthogonal group $O(p,q)$
$($$p+q: $ even, $p,q \ge 2$, and $(p,q)\ne(2,2)$$)$
has a minimal representation $\pi$ whose 
 minimal $K$-type is of the form 
$\Har{0}{p}\otimes\Har{\frac{p-q}{2}}{q}$ if $p \ge q$. 
In the $L^2$-model of $\pi$, we have proved that any $\Mmax$-fixed
vector in $\Har{0}{p}\otimes\Har{\frac{p-q}{2}}{q}$ is a scalar multiple
of the function 
$ r^{\frac{3-q}{2}}K_{\frac{q-3}{2}}(2r)$
(see \cite[Theorem 5.5]{xkors3} for a precise statement), 
where $K_\nu(z)$ denotes the $K$-Bessel function.
Since $K_{-\frac{1}{2}}(2r)= \frac{\sqrt{\pi}}{2}
r^{-\frac{1}{2}}e^{-2r}$ and $L_0^\alpha(x)=1$, we have
$ r^{\frac{3-q}{2}}K_{\frac{q-3}{2}}(2r)=\frac{\sqrt{\pi}}{2}e^{-2r}
= \frac{\sqrt{\pi}}{2} f_{0,0} (r)$
if $q=2$. 
This vector is a generator of the one dimensional vector 
space $W_{0,0}$.
\end{remNonumber} 

\begin{remNonumber}[Weil representation]\label{Remark1}
Let us compare our representation on
$L^2(\mathbb{R}^m,|x|^{-1}dx)$ with
 the (original) Schr\"{o}dinger model on $L^2(\mathbb{R}^m)$ of
 the Weil representation of $G'=\Mp$. The counterpart to Proposition 
\ref{branching-law} can be stated as follows:
we set for $0 \le l \le a$
\begin{equation}
\label{eqn:faldash}
  f'_{a,l}(r):= L_{a-l}^{\frac{m-2}2+l}(r^2) r^l e^{-\frac{r^2}2},
\end{equation}
and define linear maps by
\begin{equation*}
  j'_{a,l}: \Har{l}{m} \to C^\infty (\mathbb R^m), \quad
  \phi(\omega) \mapsto (j'_{a,l}\phi)(r\omega)=f'_{a,l}(r) \phi(\omega).
\end{equation*}
Then for any $\phi \in \Har{l}{m}$, $j'_{a,l}(\phi)$ 
is square integrable on $\mathbb R^m$, 
and its image $j'_{a,l}(\Har{l}{m})$ is characterized by the following properties:
Let $(K', R')=(U(m)\widetilde{} \; , O(m)\widetilde{} \; )$. 

1) $j'_{a,l}(\mathcal{H}^l(\mathbb{R}^m))$ 
is isomorphic to $\Har{l}{m}$ as $R'$-modules, 

2) it is contained in the $K'$-type isomorphic to $S^a(\mathbb C) \otimes 
\operatorname{det}^\frac{1}4$. 
\end{remNonumber}

The remaining part of this section is organized as follows.
In Subsection \ref{The vector Z}, we shall define a central element $Z$ of 
(the complexification of ) the Lie algebra $\mathfrak k$ of $K$, and compute its
differential action $d\pi_+(Z)$ (see Lemma \ref{differential operator d pi}).
By using this explicit form of $d\pi_+(Z)$, we prove Proposition \ref{branching-law} 
in Subsection \ref{Branching law of pi}. 
Finally, in Subsection \ref{semigroup}, by looking at the 
eigenvalues of $d\pi_+(Z)$ (see Lemma \ref{eigenvector}), we shall see that
$\set{\pi_+(e^{tZ}):=\exp (t d\pi_+(Z))}{\operatorname{Re} t >0}$ forms a 
holomorphic semigroup of contraction operators.

In Section \ref{Main Theorem}, we shall find an explicit integral kernel of this semigroup 
by using Proposition \ref{branching-law}.

\subsection{Description of infinitesimal generators of $\mathfrak{sl}(2, \mathbb R)$}
\label{The vector Z}

Let 
$$
  e:=\begin{pmatrix} 0  & 1 \\ 0  &  0 \end{pmatrix},\quad 
  f:=\begin{pmatrix} 0  & 0 \\ 1  &  0 \end{pmatrix}, \quad 
  h:=[e,f]=\begin{pmatrix} 1 & 0 \\ 0 & -1 \end{pmatrix}.
$$
be the standard basis of $\mathfrak{sl}(2, \mathbb R)$. 
With the notation (\ref{definition of E}), we define a Lie algebra
homomorphism 
\begin{equation}\label{eqn:phidef}
\phi: \mathfrak{sl}(2,\mathbb R) \to \mathfrak{so}(m+1, 2)
\end{equation}
by
\begin{equation}\label{definition of phi}
  \phi(e)=\overline{N}_{m+1}, \quad  \phi(f)= N_{m+1}, 
  \quad \phi(h)=-2E.
\end{equation}
In this subsection, we shall explicitly describe 
$d\pi_+(\phi(e)),
d\pi_+(\phi(f))$, and  $d\pi_+(\phi(h))$
as differential operators on $\mathbb{R}^m$.

\begin{lemSec}\label{Lemma3.2.1}
Let 
$E_x=\sum_{j=1}^m x_j \frac{\partial}{\partial x_j}$ and 
$\Delta=\sum_{j=1}^m \frac{\partial^2}{\partial x_j^2}$. 
Then we have:
\begin{alignat}{2} 
  d\pi_+(\phi(h))=&2E_x+m-1, \label{3.2.1-1}  \\
  d\pi_+(\phi(e))= &2\sqrt{-1}|x|,  \label{3.2.1-2}  \\
  d\pi_+(\phi(f))= &\frac{\sqrt{-1}}2 |x|\Delta. \label{3.2.1-3}
\end{alignat}
\end{lemSec}

\begin{remNonumber}
\label{rem:3.3.2}
Lemma \ref{Lemma3.2.1} corresponds to an analogous result for the 
Schr\"{o}dinger model of the Weil representation $(\varpi, L^2(\mathbb R^m))$ of $\Mp$ as follows:
by the matrix realization of the real symplectic Lie algebra
$\mathfrak{sp}(m,\mathbb{R})$, 
we define a Lie algebra homomorphism
 $\varphi: \mathfrak{sl}(2, \mathbb R) \to
\mathfrak{sp}(m, \mathbb R)$ by
\begin{equation}
  \varphi(e)=\begin{pmatrix} 0  & I_m \\ 0  &  0 \end{pmatrix},\quad 
  \varphi(f)=\begin{pmatrix} 0  & 0 \\ I_m  &  0 \end{pmatrix}, \quad 
  \varphi(h)=\begin{pmatrix} I_m & 0 \\ 0 & -I_m \end{pmatrix}.
\end{equation} 
Then, $\{ d\varpi (\varphi(h)), d\varpi(\varphi(f)),
d\varpi (\varphi (f))\}$ 
is no other than the $\mathfrak{sl}_2$-triple of differential operators
$\{ \tilde h', \tilde e', \tilde f' \}$ on $\mathbb{R}^m$ given in 
\eqref{eqn:sl2Weil}. 
\end{remNonumber}

\begin{proof}[Proof of Lemma \ref{Lemma3.2.1}]
First we compute $d\pi_+(\phi(h))=d\pi_+(-2E)$.  
For this,
we use the formula of $d\hat{\varpi}$ on $\mathbb{R}^{m+1}$ in
Subsection \ref{Infinitesimal action.},
and then compute the formula of $d\pi_+$ on the positive cone
$C_+$ (or on the coordinate space $\mathbb{R}^m$)
through the embedding
$ \imath: L^2(C) \to \mathcal{S}' (\mathbb{R}^{m+1})$,
$u(\zeta) \mapsto u(\zeta) \delta(Q)$.
Since the distribution $\delta(Q)$ is homogeneous of degree $-2$,
we note that $E_\zeta \delta(Q)=-2 \delta(Q)$.
Therefore,
\begin{alignat*}{2}
(d\pi_+(-2E)u) \delta(Q)=& -2 d\hatpi (E) (u\delta(Q)) \quad 
\text{by (\ref{commutative diagram})}
\\
=&(2E_\zeta+m+3)(u\delta(Q)) \quad \text{by (\ref{d varphi E})} \\
=&(2E_\zeta +m-1)u\cdot \delta(Q).
\end{alignat*}
Now by identifying $L^2(C_+)$ with 
$\WL$ by (\ref{jmath}), we obtain (\ref{3.2.1-1}). 

The second formula (\ref{3.2.1-2}) follows immediately from 
(\ref{d varpi N j 1}). 

We shall show the third formula (\ref{3.2.1-3}). In light of 
$\phi(f)=N_{m+1}$ (see (\ref{definition of phi})), by (\ref{d varpi N j 2}),
we have
\begin{equation}
\label{differential operator d pi 1}
d\hat{\varpi}(\phi(f)) 
=\sqrt{-1} \Bigl(\frac{m+3}2 \frac{\partial}{\partial \zeta_{m+1}} +E_\zeta 
\frac{\partial}{\partial \zeta_{m+1}}+\frac{\zeta_{m+1}}2 \square_\zeta
\Bigr).
\end{equation}
In order to compute the action $d\hat{\varpi} (\phi(f))$ along the cone $C_+$,
we use the following coordinate on $\mathbb R^{m+1}$:
\begin{equation}\label{another coordinate}
  \mathbb R \times \mathbb R_+ \times S^{m-1} \to \mathbb R^{m+1}, \quad 
  (Q, r, \omega) \mapsto (r\omega, \sqrt{r^2 -Q}). 
\end{equation}

\begin{claim}
\label{claim to differential operator d pi}
With the above coordinate, the differential operator $d\hat{\varpi}(\phi(f))$ on 
$\mathcal S'(\mathbb R^{m+1})$ takes the form: 
$$
  d\hat{\varpi}(\phi(f))=\sqrt{-1}\sqrt{r^2-Q}
  \Bigl( \frac{1}2 \frac{\partial^2}{\partial r^2}+\frac{m-1}{2r} \frac{\partial}{\partial r}
  +\frac{\Delta_{S^{m-1}}}{2r^2}- 2Q \frac{\partial^2}{\partial Q^2}  
   -4\frac{\partial}{\partial Q}\Bigr).
$$
\end{claim}
\begin{proof}[Proof of Claim \ref{claim to differential operator d pi}]
We start with a new coordinate $\mathbb R^{m+1}=\mathbb R^m\oplus
\mathbb R$ by
\begin{equation}\label{still another coordinate}
  \mathbb R_+ \times S^{m-1} \times \mathbb R \to
  \mathbb R^{m+1}, \quad (R, \omega, \zeta_{m+1}) \mapsto (R\omega, \zeta_{m+1}).
\end{equation}
Then, clearly,
\begin{alignat*}{1}
  E_\zeta =&R\frac{\partial}{\partial R}+
\zeta_{m+1} \frac{\partial}{\partial \zeta_{m+1}} \\
\square_\zeta = &\Delta_{\mathbb R^m} -\frac{\partial^2}{\partial \zeta_{m+1}^2} 
                           =\Bigl( \frac{\partial^2}{\partial R^2}
                                        +\frac{m-1}R \frac{\partial}{\partial R} 
                                        +\frac{1}{R^2} \Delta_{S^{m-1}} \Bigr)
                              -\frac{\partial^2}{\partial \zeta_{m+1}^2}.    
\end{alignat*}
The coordinate (\ref{another coordinate}) is obtained by the composition of 
(\ref{still another coordinate}) and
$$
\quad r=R, \quad Q=R^2-\zeta_{m+1}^2.
$$
In light of
\begin{equation}
\label{claim 1}
\frac{\partial}{\partial R}=\frac{\partial}{\partial r}+2r \frac{\partial}{\partial Q},
\quad 
\frac{\partial}{\partial \zeta_{m+1}}=-2\sqrt{r^2-Q} \frac{\partial}{\partial Q},
\end{equation}
we get
\begin{alignat}{1}
\label{claim 2}
     E_\zeta =& r\frac{\partial}{\partial r}+ 2Q\frac{\partial}{\partial Q},  \\
\label{claim 3}
      \square_\zeta=& \frac{\partial^2}{\partial r^2} 
                                  +4Q \frac{\partial^2}{\partial Q^2}
                                  +4r\frac{\partial^2}{\partial r \partial Q}
                                  +2(m+1)\frac{\partial}{\partial Q}
                                  +\frac{m-1}{r}\frac{\partial}{\partial r}
                                  +\frac{\Delta_{S^{m-1}}}{r^2}.
\end{alignat}
Now substituting (\ref{claim 1}), (\ref{claim 2}) and (\ref{claim 3}) into 
(\ref{differential operator d pi 1}), we get the claim.
\end{proof}

Given $u \in L^2(C_+)_K$,
we extend it to a distribution 
$\tilde{u} \equiv \tilde{u}(Q,r,\omega) \in
 \mathcal{S}'(\mathbb{R}^{m+1}) \cap C^\infty(\mathbb{R}^m \setminus
  \{0\})$ 
such that
$$
u \delta(Q) = \tilde{u} \delta(Q).
$$
We set
$$
  S:=2 \sqrt{r^2-Q} \Bigl( Q \frac{\partial^2}{\partial Q^2} +
    2\frac{\partial}{\partial Q}\Bigr)
        \bigl(\Tilde u \delta (Q)\bigr). 
$$
Then, it follows from \eqref{commutative diagram} and 
 Claim \ref{claim to differential operator d pi} that 
\begin{alignat*}{1}
&d\hat{\varpi}(\phi(f))(\Tilde u \delta(Q))
\\
&=
\sqrt{-1}\Bigl( \sqrt{r^2-Q} \Bigl(\frac{1}2 \frac{\partial^2}{\partial r^2} +\frac{m-1}{2r}
 \frac{\partial}{\partial r}+\frac{\Delta_{S^{m-1}}}{2r^2}\Bigr)\Tilde u \Bigr) \delta(Q) 
-S   \\
&=\sqrt{-1}\Bigl( 1+O(\bigl(\frac{Q}{r^2} \bigr) ) \Bigr)
     \Bigl(\frac{r}{2}\frac{\partial^2}{\partial r^2}  
               +\frac{m-1}{2} \frac{\partial}{\partial r}
               +\frac{\Delta_{S^{m-1}}}{2r}\Bigr) u\cdot \delta(Q)-S \\
&=\sqrt{-1}\Bigl(\frac{r}{2}\frac{\partial^2}{\partial r^2}  
               +\frac{m-1}{2} \frac{\partial}{\partial r}
               +\frac{\Delta_{S^{m-1}}}{2r}\Bigr)u \cdot \delta(Q)-S
\\
&=\frac{\sqrt{-1}}{2} r (\Delta_{\mathbb{R}^m} u) \delta(Q)-S.
\end{alignat*}
At the second last equality, we used the formula $Q \delta(Q)=0$.
In the following claim, we shall show $S=0$.
Now, 
the proof of \eqref{3.2.1-3} is complete.
Hence, we have shown Lemma \ref{Lemma3.2.1}. 
\end{proof}
\begin{claim}
\label{claim to differential operator d pi 1}
For any $\Tilde u \in C_0^\infty(\mathbb R^{m+1}\setminus \{0\})$, we have
$$
  \Bigl(Q\frac{\partial^2}{\partial Q^2}+2\frac{\partial}{\partial Q}\Bigr)
  \bigl(\Tilde u \delta(Q)\bigr)=0.
$$
\end{claim}

\begin{proof}[Proof of Claim \ref{claim to differential operator d pi 1}]
By the Leibniz rule, the left-hand side amounts to
$$
  \frac{\partial^2 \Tilde u}{\partial Q^2}
  Q \delta (Q) +2 \frac{\partial \Tilde u}{\partial Q} \Bigl( Q \frac{\partial}{\partial Q}
  \delta (Q) +\delta(Q)
  \Bigr) + \Tilde u \Bigl( Q \frac{\partial^2}{\partial Q^2} \delta (Q)+
  2\frac{\partial}{\partial Q} \delta(Q) \Bigr).
$$
Hence we see that this equals $0$ in light of the formulas:  
$$
  Q \delta(Q)=0,\quad  
  Q\frac{\partial}{\partial Q}\delta(Q)=-\delta (Q), \quad 
  Q \frac{\partial^2}{\partial Q^2} \delta(Q)= -2 \frac{\partial}{\partial Q} \delta(Q).
$$
\end{proof}

\subsection{Central element $Z$ of $\mathfrak{k}_{\mathbb{C}}$}
\label{sec:3.4}
We extend the Lie algebra homomorphism
$\phi: \mathfrak{sl}(2,\mathbb{R}) \to \mathfrak{so}(m+1,2)$
(see (\ref{definition of phi})) 
to the complex Lie algebra homomorphism
$\phi: \mathfrak{sl}(2,\mathbb{C}) \to \mathfrak{so}(m+1,\mathbb{C})$.
Consider a generator of $\mathfrak{so}(2,\mathbb{C})$ given by
\begin{equation}\label{definition of z}
  z:= \frac{\sqrt{-1}}2 (e-f)=\frac{\sqrt{-1}}2 \begin{pmatrix} 0 & 1
  \\ -1 & 0\end{pmatrix}. 
\end{equation}
We set
\begin{equation}\label{definition of Z}
Z:= \phi(z) \in \sqrt{-1} \mathfrak{g}. 
\end{equation}
In light of (\ref{definition of phi}) and (\ref{definition of E}), we have
\begin{align} \label{3.4.1}
Z&= \frac{\sqrt{-1}}2 (\overline{N}_{m+1}-N_{m+1})
\\
&=\sqrt{-1}( E_{m+1,m+2}-E_{m+2,m+1}). 
\label{eqn:ZcenterK}
\end{align}
Hence, $\sqrt{-1}Z$ is contained in the center 
$\mathfrak{c}(\mathfrak{k})$ of 
$\mathfrak{k} \simeq \mathfrak{so}(m+1) \oplus \mathfrak{so}(2)$. 
If $m>1$, then $\mathfrak{c}(\mathfrak{k})$ is of one dimension,
and $\sqrt{-1} Z$ generates $\mathfrak{c}(\mathfrak{k})$. 
By \eqref{definition of z}, we have 
$e^{\sqrt{-1}tZ}=I_{m+3}$ in $SO_0(m+1,2)$ if and only if 
$t \in 2\pi \mathbb Z$. Hence,
$e^{\sqrt{-1}tZ}=1$ in $G=SO_0(m+1,2)\widetilde{}$ if and only if
$t \in 4\pi \mathbb Z$ (see \eqref{eqn:Gdouble}).
Therefore, we have the following lemma:
\begin{lemSec}\label{lem:SL2G}
The Lie algebra homomorphism 
$\phi: \mathfrak{sl}(2,\mathbb R)\to \mathfrak{so}(m+1,2)$
(see \eqref{eqn:phidef}) lifts to an injective Lie group 
homomorphism $SL(2,\mathbb R)\to G$.
\end{lemSec}
By the expression \eqref{eqn:ZcenterK} of $Z$ and by the
$K$-isomorphism \eqref{eqn:Wa},
we have:

\begin{lemSec}
\label{lem:infZaction}
$d\pi_+(Z)$ acts on $W_a$ as a scalar multiplication of\/ 
$- (a+\frac{m-1}{2} )$. 
\end{lemSec}

Combining (\ref{definition of z}), (\ref{definition of Z}), and Lemma \ref{Lemma3.2.1}, 
we readily get the following lemma:
\begin{lemSec}
\label{differential operator d pi}
The differential operator $d\pi_+(Z)$ takes the form:
$$
  d\pi_+(Z)=|x| \Bigl( \frac{\Delta}4-1 \Bigr).
$$
\end{lemSec}

\begin{remNonumber}\label{Remark}
$d\pi_+(Z)$ coincides with the operator $D$ in Introduction. 
In particular, $D$ is a self-adjoint operator on $\WL$ because $\sqrt{-1}Z 
\in \mathfrak g$ and $\pi_+$ is a unitary representation.
\end{remNonumber}

\begin{remNonumber}[Weil representation]\label{Remark2}
For the Weil representation of $G'$, 
the Lie algebra homomorphism
$\varphi: \mathfrak{sl}(2,\mathbb{R}) \to \mathfrak{sp}(m,\mathbb{R})$
(see Remark \ref{rem:3.3.2}) sends $z$ to
$$
  Z':=\varphi(z) = \frac{\sqrt{-1}}2 \begin{pmatrix}  0 &  I_m  \\ -I_m  & 0 \end{pmatrix}
 \in \sqrt{-1} \mathfrak g'.
$$ 
Hence $\sqrt{-1}Z'$ is a central element of $\mathfrak k' \simeq 
\mathfrak u(m)$.
The differential operator $d\varpi(Z')$ amounts to the Hermite operator
(see \cite[\S6 (d)]{xHo})
$$
\mathcal D=\frac{1}4 (\Delta -|x|^2).
$$ 
\end{remNonumber}

See the table in 
Subsection \ref{Minimal representation as hidden symmetry} 
for the differential operators corresponding to $e,f$ and $h$.

\subsection{Proof of Proposition \protect\ref{branching-law}}
\label{Branching law of pi}

This subsection gives a proof of Proposition \ref{branching-law}.

1) Let us show 
$j_{a,l}(\phi)=f_{a,l}\phi \in 
L^1(\mathbb R^m, \frac{dx}{|x|})\cap
\WL$ 
for any $\phi \in \Har{l}{m}$. 
By the definition (\ref{definition of f_a,l}) of $f_{a,l}$, 
$f_{a,l}(r)$ is regular at $r=0$, 
and $f_{a,l}(r)$ decays exponentially as $r$ tends to infinity. 
Therefore, $f_{a,l} \in 
L^1(\mathbb R_+, r^{m-2}dr) \cap
\wL$. Since our measure $\frac{dx}{|x|}$ has the form 
\begin{alignat}{2}
\label{4.1.0}
  \frac{dx}{|x|}=r^{m-2}dr d\omega, 
\end{alignat}
with respect to the polar coordinate (\ref{polar coordinate}), we have
shown 
$j_{a,l}(\phi) \in L^1(\mathbb R^m, \frac{dx}{|x|}) \cap \WL$.

2) We set $H_{a,l}:=j_{a,l}\bigl(\Har{l}{m} \bigr)$. 
Obviously, $H_{a,l}$ is isomorphic to $\mathcal{H}^l(\mathbb{R}^m)$ as
$R$-modules. 
To see $W_{a,l} = H_{a,l}$,
it is sufficient to show the following inclusion:
\begin{equation}
\label{inclusion}
  H_{a,l} \subset W_a,
\end{equation}
because $W_{a,l}$ is characterized as the unique 
subspace of $W_a$ 
such that $W_{a,l} \simeq \Har{l}{m}$ as $\R$-modules. 

To see \eqref{inclusion}, 
we recall from \eqref{eqn:Wa} that $W_a$ is characterized as the
unique subspace of $L^2(\mathbb{R}^m,\frac{dx}{|x|})$ 
on which $d\pi_+(Z)$ acts with eigenvalue 
$-(a+\frac{m-1}2)$. Thus the inclusive relation 
(\ref{inclusion}) will be proved if we show the following lemma:

\begin{lemSec}
\label{eigenvector}
The operator $d\pi_+(Z)$ acts on $H_{a,l}$ as a scalar multiplication $-(a+\frac{m-1}2)$. 
In other words, we have
\begin{equation}
\label{eigenvector-1}
 d\pi_+(Z)(f_{a,l}\phi) = -\Bigl(a+\frac{m-1}{2}\Bigr) f_{a,l}\phi.
\end{equation}
\end{lemSec}

\begin{proof}[Proof of Lemma \ref{eigenvector}]
Writing the differential operator $d\pi_+(Z)$ (see Lemma 
\ref{differential operator d pi}) in terms of the polar coordinate 
and using the definition of $f_{a,l}$ (see (\ref{definition of f_a,l})), we see
that the equation (\ref{eigenvector-1}) is equivalent to 
\begin{alignat}{2}
\label{eigenvector-2}
\biggl(\frac{r}{4}\frac{\partial^2}{\partial r^2}+
\frac{m-1}{4} \frac{\partial}{\partial r}
+\frac{\Delta_{S^{m-1}}}{4r} -r +\Bigl(a+\frac{m-1}2\Bigr) \biggr)
\bigl( \psi(4r)r^l e^{-2r}\phi(\omega)\bigr)  
= 0,
\end{alignat}
for $\psi(r):=L_{a-l}^{m-2+2l}(r)$.
The equation (\ref{eigenvector-2}) amounts to 
\begin{equation*}
\Bigl( 4r \psi''(4r) +(m-1+2l-4r)\psi'(4r)+(a-l)\psi(4r) \Bigr)
r^l e^{-2r} \phi(\omega)=0.
\end{equation*}
This is nothing but  Laguerre's differential equation (\ref{ode})
with $n=a-l, \alpha = m-2+2l$. Now the lemma follows.
\end{proof}

\subsection{One parameter holomorphic semigroup $\pi_+(e^{tZ})$}
\label{semigroup}
It follows from Lemma \ref{lem:infZaction} that 
for $t \in \mathbb C$ the operator
\begin{equation}
  \pi_+(e^{tZ}):= \exp d\pi_+(tZ) = \sum_{j=0}^\infty \frac{1}{j!} d\pi_+(tZ)^j
\end{equation}
acts on $W_{a,l}$ as a scalar multiplication of 
$e^{-(a+\frac{m-1}{2})t}$ 
for any $0 \le l \le a$, 
 of which the absolute value 
does not exceed $1$ if $\operatorname{Re}t \ge 0$. 
In light of the direct sum decomposition \eqref{R_0 decomposition}, 
if $\operatorname{Re} t \ge 0$ then 
 the linear map
$\pi_+(e^{tZ}): \WL_{K} \to \WL_{K}$ extends to a 
continuous operator (we use the same notation $\pi_+(e^{tZ})$) on 
$\WL$.
Furthermore, it is a contraction
operator if $\operatorname{Re} t>0$. We summarize some of basic 
properties of $\pi_+(e^{tZ})$:

\begin{propSec}
\label{Lemma3.4}
{\rm 1)} The map 
\begin{equation}
  \set{ t\in \mathbb{C} }{\operatorname{Re}t\geq 0 } \times
  \WL  \to \WL, \quad (t,f) \mapsto \pi_+(e^{tZ})f
\end{equation}
is continuous. 

{\rm 2)} For a fixed  $t$ such that $\operatorname{Re}t\geq 0$,  
$\pi_+(e^{tZ})$ is characterized as 
the continuous operator from
$\WL \to \WL$ satisfying 
\begin{equation}
\label{eigenfunction_of_T}
 \pi_+(e^{tZ}) u= e^{-(a+\frac{m-1}{2})t}u, 
\end{equation}
for any 
$u \in W_{a,l}=\{f_{a,l}\phi: \phi\in\mathcal{H}^l(\mathbb{R}^m)\}$ 
(see Proposition \ref{branching-law}) and 
for any $l, a \in \mathbb N$ such that $0 \le l \le a$.

{\rm 3)} The operator norm $\| \pi_+(e^{tZ})\|$
of $\pi_+(e^{tZ})$ is $e^{-\frac{m-1}2 \operatorname{Re}t}$. 

{\rm 4)} If $\operatorname{Re}t>0$, $\pi_+(e^{tZ})$ is a Hilbert--Schmidt operator.

{\rm 5)} If $t \in \sqrt{-1} \mathbb R$, then 
$\pi_+(e^{tZ})$ is a unitary operator.
\end{propSec}

\begin{remNonumber}
We define a subset $\Gamma_+:=\set{ tZ}{t>0}$ in $\sqrt{-1}\mathfrak g$.
Proposition \ref{Lemma3.4} indicates how the unitary representation $\pi_+$ of $G$
extends to a holomorphic semigroup on the 
complex domain $G \cdot \exp \Gamma_+ \cdot G$ 
of $G_\mathbb C$. Our results may be regarded as a part of 
the {\it Gelfand--Gindikin program}, which tries to understand 
unitary representations of a real semisimple Lie group
by means of holomorphic objects on an open subset of $G_\mathbb C \setminus G$,
where $G_\mathbb C$ is a complexification of $G$ 
(see \cite{xGeGi,xOl,xSt}).
\end{remNonumber}

\begin{remNonumber}[Weil representation]
In the case of the Weil representation $(\varpi, L^2(\mathbb R^m))$ of 
$G'=\Mp$,
the functions $f'_{a,l}$ (see Remark \ref{Remark1}) play the same role as 
$f_{a,l}$ because of the following facts:

1) $\set{f'_{a,l}}{0\le l \le a}$ spans a complete orthogonal
   basis
 of $L^2(\mathbb R_+, r^{m-1}dr)$ (cf.\ Remark \ref{rem:falcob}).

2) For any $\phi \in \Har{l}{m}$, 
$f'_{a,l}(r)\phi(\omega)$ $(0 \le l \le a)$ are eigenfunctions of $d\varpi(Z')$ 
with negative eigenvalues $-(\frac{m}4+ \frac{l}2+\frac{a}2)$.

Owing to these facts, we obtain a holomorphic semigroup of contraction operators
$\varpi(e^{tZ'}):=\exp t(d\varpi(Z'))$ $(\operatorname{Re}t>0)$ on $L^2(\mathbb R^m)$.
Since $d\varpi(Z')$ coincides with the Hermite operator $\mathcal D$
(see Remark \ref{Remark2}), this holomorphic semigroup is 
nothing but the Hermite semigroup (see \cite[\S5]{xHo}).
\end{remNonumber}

\section{Radial part of the semigroup}
\label{Radial part of the semigroup}

This section gives an explicit integral formula for the \lq radial part\rq\ 
of the holomorphic semigroup $\pi_+(e^{tZ})$ in the `Schr\"{o}dinger
model' $L^2(\mathbb{R}^m, |x|^{-1} dx)$.
The main result of Section \ref{Radial part of the semigroup} is Theorem 
\ref{Theorem-C}. 
As its applications, we see that the semigroup
law 
$ \pi_+(e^{(t_1+t_2)Z})=\pi_+(e^{t_1Z})\circ \pi(e^{t_2Z})$
gives a
simple and representation theoretic proof of the classical 
Weber's second exponential integral 
formula on Bessel functions (Corollary \ref{Corollary-B}),
and that taking the boundary value 
$\lim\limits_{s\downarrow 0} \pi_+(e^{sZ}) = \operatorname{id}$
provides an example of
 a Dirac sequence 
 (Corollary \ref{Corollary-Dirac}). 

Theorem \ref{Theorem-C} will play a key role in 
Section \ref{Main Theorem},
where we complete the proof of 
the main theorem of this article, namely, 
Theorem \ref{Theorem-B} that
gives an integral formula of the holomorphic semigroup
$\pi_+(e^{tZ})$ on $L^2(\mathbb{R}^m, |x|^{-1} dx)$.

\subsection{Result of the section}
\label{Result of Radial part of the semigroup}
For a complex parameter $t \in \mathbb C$ with $\operatorname{Re}t>0$,
we have defined a contraction operator 
$\pi_+(e^{tZ}):\WL \to \WL$ in Proposition \ref{Lemma3.4}.

We recall that
 $Z$ is a central element in $\mathfrak{k}_{\mathbb{C}}$
(see Subsection \ref{sec:3.4}) 
and that $R$ is a subgroup of $K$.
Therefore, 
$\pi_+(e^{tZ})$ intertwines with the $\R$-action.
On the other hand, the natural action of
$R \simeq SO(m)$ gives a direct sum decomposition of the Hilbert
space: 
\begin{equation} \label{4.1.1}
  \WL \simeq 
 \sideset{}{^\oplus}{\sum}_{l=0}^{\infty}
  \wL \otimes \Har{l}{m}. 
\end{equation}
Hence, by Schur's lemma, 
there exists a family of continuous operators parametrized
by $l \in \mathbb{N}$:
\begin{equation} \label{definition of T_l}
  \pi_{+, l}(e^{tZ}): \wL \to \wL 
\end{equation}
such that  $\pi_+(e^{tZ})$ is diagonalized
according to the direct sum decomposition \eqref{4.1.1} as follows:
\begin{equation}\label{eqn:pilsum}
  \pi_+(e^{tZ})= \sideset{}{^\oplus}{\sum}_{l=0}^\infty
  \pi_{+,l}(e^{tZ})\otimes \operatorname{id}. 
\end{equation}

The goal of this section is to give an explicit integral formula
of  $\pi_{+,l}(e^{tZ})$ on $\wL$ for $l\in\mathbb{N}$.
We note that $\pi_{+,l}(e^{tZ})$ is a unitary 
operator if $\operatorname{Re} t=0$ because so is
$\pi_+(e^{tZ})$. 
Likewise, $\pi_{+,l}(e^{tZ})$ is a Hilbert--Schmidt operator if
$\operatorname{Re} t> 0$ because so is $\pi_+(e^{tZ})$. 

We now introduce the following subset of $\mathbb C$: 
\begin{equation}
\label{definition of D}
  \Omega:=\set{t \in \mathbb C}{\operatorname{Re}t \ge 0}
        \setminus 2\pi\sqrt{-1} \mathbb Z, 
\end{equation}
and define a family of analytic functions $K_l^+(r, r'; t)$ on $\mathbb R_+
\times \mathbb R_+ \times \Omega$ by the formula:
for $l=0,1,2,\ldots\,$,
\begin{alignat}{2}
  K_l^+(r,r';t):={}&\frac{2 e^{-2(r+r')\coth\frac{t}{2}}    }{\sinh\frac{t}{2}}
          (rr')^{-\frac{m-2}{2}}I_{m-2+2l}
          \Bigl(\frac{4\sqrt{rr'} }{\sinh\frac{t}{2}} \Bigr)
\nonumber
\\
={}&\frac{2^{m-1+2l} e^{-2(r+r')\coth\frac{t}{2}} (rr')^l}
         {(\sinh\frac{t}{2})^{m-1+2l}}
    \tilde{I}_{m-2+2l}
    \Bigl(\frac{4\sqrt{rr'} }{\sinh\frac{t}{2}} \Bigr).
\label{equation3}
\end{alignat}
Here, $\tilde{I}_\nu(z) := (\frac{z}{2})^{-\nu} I_\nu(z)$, and
$I_\nu(z)$ denotes the $I$-Bessel function (see Subsection \ref{Bessel functions}).
We note that the denominator $\sinh \frac{t}2$ is nonzero
for $t \in \Omega$.

We are ready to state the
 integral formula
of  $\pi_{+,l}(e^{tZ})$ for $t \in \Omega$:
\begin{thmSec}[Radial part of the semigroup]
\label{Theorem-C}
{\rm 1)} For $\operatorname{Re}t>0$, 
the Hilbert--Schmidt operator $\pi_{+,l}(e^{tZ})$ 
 on $\wL$
 is given by the following integral transform: 
\begin{equation}
\label{4.2.1}
  (\pi_{+,l}(e^{tZ})f)(r)= \int_0^\infty K_l^+(r, r'; t)f(r') {r'}^{m-2}dr'.
\end{equation}
The right-hand side converges absolutely
for $f \in \wL$.

{\rm 2)}
If $t \in \sqrt{-1}\mathbb R$ but $t\notin 2\pi\sqrt{-1}\mathbb{Z}$, 
then the integral formula \eqref{4.2.1} for the unitary operator
$\pi_{+,l}(e^{tZ})$ 
holds in the sense of $L^2$-convergence. 
Furthermore, the right-hand side converges absolutely 
if $f$ is 
a finite linear combination 
of $f_{a,l}$ $(a=l, l+1, \cdots)$. 
\end{thmSec}

\begin{remNonumber}
Let us compare Theorem \ref{Theorem-C} with the corresponding result
for the Weil representation $\varpi$
of $G'=\Mp$ realized as the Schr\"{o}dinger model $L^2(\mathbb{R}^m)$. 
According to the direct sum decomposition of the Hilbert space:
$$
  L^2(\mathbb R^m) \simeq \sideset{}{^\oplus}{\sum}_{l=0}^\infty
  L^2(\mathbb R_+, r^{m-1}dr) \otimes \Har{l}{m},
$$
there exists a family of continuous operators
$\varpi_l(e^{tZ'})$ $(\operatorname{Re} t>0)$ such that 
the holomorphic semigroup $\varpi(e^{tZ'})$ 
has the following decomposition: 
$$
  \varpi (e^{tZ'}) = \sideset{}{^\oplus}{\sum}_{l=0}^\infty
  \varpi_l(e^{tZ'}) \otimes \operatorname{id}.
$$
Then, by an analogous computation 
to Theorem \ref{Theorem-C}, 
we find the kernel function of the 
semigroup $\varpi_l(e^{tZ'})$ is given by 
$$
  \mathcal K_l(r, r'; t):= \frac{e^{-\frac{1}2(r^2+{r'}^2)\coth \frac{t}2}}{\sinh \frac{t}2}
                        (rr')^{-\frac{m-2}2} I_\frac{m-2+2l}2\Bigl(\frac{rr'}{\sinh \frac{t}2}\Bigr).
$$

The relation between the kernel function $\mathcal{K}$
(Mehler kernel) of $\varpi(e^{tZ'})$ and
$\mathcal{K}_l(r, r'; t)$ will be discussed in Remark \ref{rem:5.7.2}.
\end{remNonumber}

This section is organized as follows. 
In Subsection \ref{sec:4.2}, we give a proof of Theorem \ref{Theorem-C} for the 
case $\operatorname{Re}t>0$, which is based on a computation of the kernel 
function by means of the infinite sum of the eigenfunctions. 
In Subsection \ref{sec:4.3}, by taking the analytic continuation, the case 
$\operatorname{Re}t=0$ is proved. 
Applications of Theorem \ref{Theorem-C} to special function theory are
discussed in Subsections \ref{sec:4.4} and \ref{Dirac sequence operators}.

\subsection{Upper estimate of the kernel function}

In this subsection, we shall
give an upper estimate of
the kernel function $K_l^+(r,r';t)$.

For $t = x+\sqrt{-1}y$, we set
\begin{align}
&\alpha(t) := \frac{\sinh x}{\cosh x - \cos y},
\label{eqn:atxy}
\\
&\beta(t) := \frac{\cos\frac{y}{2}}{\cosh\frac{x}{2}}.
\label{eqn:btxy}
\end{align}
Then, an elementary computation shows
\begin{align}
&\operatorname{Re} \coth\frac{t}{2} = \alpha(t),
\label{eqn:Reat}
\\
&\operatorname{Re} \frac{1}{\sinh\frac{t}{2}}
  = \alpha(t) \beta(t).
\label{eqn:Reabt}
\end{align}
For $t \in \Omega$
(see \eqref{definition of D} for definition),
we have $\cosh x - \cos y > 0$, and then,
\begin{equation}\label{eqn:abpos}
  \alpha(t)\ge 0 \quad \text{and}\quad |\beta(t)|<1.
\end{equation}
If $\operatorname{Re}t>0$, then
\begin{equation}
  \alpha(t)>0.
\end{equation}
For later purposes, we prepare:

\begin{lemSec}\label{lem:estimate}
If $y\in \mathbb R$ satisfies $|y|\le 4\sqrt{rr'}$,
then, for $t\in \Omega$, we have the following estimate
for some constant $C$:
\begin{equation}
  \bigl|e^{-2(r+r')\coth \frac{t}2}\tilI_\nu (\frac{y}{\sinh \frac{t}2})
  \bigr| \le C e^{-2\alpha(t)(1-|\beta(t)|)(r+r')}.
\end{equation}

\end{lemSec}
\begin{proof}
Using the upper estimate of the $I$-Bessel function (see Lemma
  \ref{4.3.1}),
$$
  |\tilI_\nu(\frac{y}{\sinh \frac{t}2})|\le
  C e^{|y||\operatorname{Re}\frac{1}{\sinh\frac{t}2}|},
$$
we have
\begin{alignat*}{1}
   \bigl|e^{-2(r+r')\coth \frac{t}2}\tilI_\nu (\frac{y}{\sinh \frac{t}2})
  \bigr|\le &
  C e^{-2(r+r')\operatorname{Re}\coth \frac{t}2+|y|
  \operatorname{Re}\bigl|\frac{1}{\sinh \frac{t}2}\bigr|} \\
  \le & C e^{-2(r+r')\alpha(t)+4\sqrt{rr'}\alpha(t)|\beta(t)|} \\
  \le & C e^{-2\alpha(t)(1-|\beta(t)|)(r+r')}.
\end{alignat*}
Here, the last inequality follows from
\begin{equation}
\label{eqn:muladd}
r+r'-2|\beta(t)|\sqrt{rr'}
  \ge (1-|\beta(t)|)(r+r')
\end{equation}
for $t \in \Omega$. Thus Lemma is proved.
\end{proof}
Now we state a main result of this subsection:
\begin{lemSec}\label{lem:Klupper}
Let $l\in \mathbb N$ and $m\ge 2$.

{\rm 1)} There exists a constant $C>0$ such that
\begin{equation}\label{eqn:Klupper}
  |K_l^+(r,r';t)|\le \frac{C(rr')^l e^{-2\alpha(t)(1-|\beta(t)|)(r+r')}}
           {|\sinh\frac{t}{2}|^{m-1+2l}}.
\end{equation}
for any $r, r'\in \mathbb R_+$ and $t\in \Omega$.

{\rm 2)} If $\operatorname{Re}t>0$, then
$K_l^+(\cdot,\cdot;t)\in L^2((\mathbb R_+)^2, (rr')^{m-2}drdr')$.

{\rm 3)} If $\operatorname{Re}t>0$, then for a fixed $r>0$, we have
$K_l^+(r,\cdot;t)\in L^2(\mathbb R_+, r'^{m-2}dr')$.
\end{lemSec}
\begin{proof}
1) By the definition of $K_l^+$ (see \eqref{equation3}), we have
\begin{align}\label{eqn:Klabt}
  |K_l^+(r,r';t)|=& \frac{2^{m-1+2l}(rr')^l}{|\sinh \frac{t}2|^{m-1+2l}}
   \bigl|e^{-2(r+r')\coth \frac{t}2}\tilI_{m-2+2l}
  (\frac{4\sqrt{rr'}}{\sinh \frac{t}2})
  \bigr|. \nonumber
\end{align}
Now \eqref{eqn:Klupper} follows from Lemma \ref{lem:estimate}
by substituting $\nu=m-2+2l$ and $y=4\sqrt{rr'}$.

Since $\alpha(t)>0$ for $\operatorname{Re}t>0$, the statements
2) and 3) hold by 1).
\end{proof}

\subsection{Proof of Theorem \protect\ref{Theorem-C} (Case $\operatorname{Re}t>0$)}
\label{sec:4.2}

We recall from Remark \ref{rem:falcob} that
$$f_{a,l}(r) = L_{a-l}^{m-2+2l} (4r)r^l e^{-2r}
\quad (a = l, l+1, \ldots)
$$
forms a complete orthogonal basis of
$L^2(\mathbb{R}_+, r^{m-2} dr)$.
Further, by the orthogonal relation of the Laguerre polynomials
$L_m^\alpha(x)$
(see \eqref{eqn:Lagnorm}),
we have the normalization of $\{ f_{a,l} \}$ as follows:

\begin{lemSec}\label{Lemma3.1.1}
For integers $a,b \ge l$, we have
\begin{equation}\label{3.1.1-1}
  \int_0^\infty f_{a,l}(r) f_{b,l}(r) r^{m-2}dr=
  \begin{cases} 
    0 \quad &\text{if $a\neq b$,} \\
   \frac{\Gamma(m-1+a+l)}{4^{m-1+2l} \Gamma(a-l+1)} 
   \quad &\text{if $a=b$}.
  \end{cases}
\end{equation}
\end{lemSec}

We rewrite (\ref{R_0 decomposition}) by using Proposition \ref{branching-law}
as follows:
\begin{align} 
  \WL_{K} &=\bigoplus_{a=0}^\infty\bigl(\bigoplus_{l=0}^a  W_{a,l} \bigr)
\nonumber
\\
          &=\bigoplus_{l=0}^\infty \bigl( \bigoplus_{a=l}^\infty 
                                               W_{a,l} \bigr)
\nonumber
\\
          &=\bigoplus_{l=0}^\infty \Bigl( \bigl( \bigoplus_{a=l}^\infty
            \mathbb C f_{a,l} \bigr) \otimes \Har{l}{m} \Bigr) .
\label{4.4}
\end{align}
It follows from Proposition \ref{Lemma3.4} (2) and the 
definition \eqref{eqn:pilsum} of $\pi_{+,l}(e^{tZ})$ that 
\begin{equation*}
  \pi_{+,l}(e^{tZ})f_{a,l}=e^{-(a+\frac{m-1}2)t} f_{a,l}  \quad (a=l,l+1, \cdots).
\end{equation*}

Therefore, the kernel 
function 
$K_l^+(r,r';t)$
of $\pi_{+,l}(e^{tZ})$ can be written as the infinite sum:
\begin{alignat}{1} \label{4.5}
  K_l^+(r, r'; t) = &\sum_{a=l}^\infty 
  \frac{ e^{-(a+\frac{m-1}2)t} f_{a,l}(r) \overline{f_{a,l}(r')}}
           { || f_{a,l} ||^2_{\wL} }. 
\end{alignat}
Since $\pi_{+,l}(e^{tZ})$ is a Hilbert--Schmidt operator if $\operatorname{Re}t>0$, 
the right-hand side converges in 
$L^2((\mathbb{R}_+)^2, (rr')^{m-2}drdr')$,
and therefore converges for almost all
$(r,r') \in (\mathbb{R}_+)^2$.

Let us compute the infinite sum (\ref{4.5}). 
For this,
we set
\begin{equation*}
  \kappa(r, r'; t):=\sum_{a=l}^\infty \frac{\Gamma(a-l+1)}{\Gamma(m-1+a+l)}
        L_{a-l}^{m-2+2l}(4r) L_{a-l}^{m-2+2l}(4r') e^{-(a-l)t}  .
\end{equation*}
Then, it follows from
\eqref{3.1.1-1} that we have
\begin{equation}
\label{4.6}
K_l^+(r,r';t) = 
     4^{m-1+2l} (rr')^l e^{-2(r+r')} e^{-(l+\frac{m-1}2)t} \kappa(r, r'; t).
\end{equation}
Now, we apply the Hille--Hardy formula (see \eqref{cor-A}) with 
$\alpha=m-2+2l, n=a-l, x=4r, y=4r'$, and $w=e^{-t}$.
We note that
$|w|=e^{-\operatorname{Re}t}<1$ by the assumption $\operatorname{Re}t>0$.
Then we have 
\begin{alignat*}{1}
\kappa(r,r';t)
  ={}&
   \frac{e^{\frac{(4r+4r')e^{-t}}{1-e^{-t}}}
         (-16 rr' e^{-t} )^{-\frac{m-2}2+l} }
        {1-e^{-t}}
  J_{m-2+2l} \Bigl(\frac{2\sqrt{-16 rr'e^{-t}}}{1-e^{-t}} \Bigr) \\
  ={}& \frac{(rr')^{-\frac{m-1}2+l} e^{-2(r+r')\frac{2e^{-t}}{1-e^{-t}}} e^{(\frac{m-2}2+l)t}}
                {4^{m-2+2l} (1-e^{-t})}
        I_{m-2+2l} \Bigl( \frac{4\sqrt{rr'}}{\sinh\frac{t}2} \Bigr).
\end{alignat*}
Hence, the formula (\ref{equation3}) is proved.

Therefore, the right-hand side of (\ref{4.2.1}) converges absolutely 
 by the Cauchy--Schwarz inequality because
$K_l^+(r,r';t) \in L^2(\mathbb R_+, {r'}^{m-2}dr')$
for any $r>0$ and $\operatorname{Re}t > 0$
(see Lemma \ref{lem:Klupper} (3)).

\begin{remNonumber}
\label{rem:4.2.2}
The special functions and related formulas that arise in the analysis of 
the radial part have a scheme of generalization from
$\SL$ to $G=SO_0(m+1,2)\widetilde{}$ and
$G'=\Mp$. 
This scheme is illustrated as follows:
$$
\begin{matrix}
  \SL   &  \Rightarrow  &  G ~\text{or}~G'                                                                                                            
  \\
  \text{ Segal--Shale--Weil representation } & \Rightarrow & \text{ minimal representation}
  \\
  \text{ Hermite polynomials}
  & \Rightarrow &  \text{ Laguerre polynomials}
  \\
  \text{  Mehler's formula }  & \Rightarrow  &  \text{ Hille--Hardy formula}  
\end{matrix}
$$
See \cite[p 116, Exercise 5 (d)]{xHoTa} for the $SL(2, \mathbb R)$
case.
Owing to the reduction formula of Laguerre polynomials to Hermite
polynomials (see \eqref{eqn:LHeven} and \eqref{eqn:LHodd}),
the radial part $f'_{a,l}(r)$
(see Remark \ref{Remark1}) for the Weil representation of
$G' = \Mp$ collapses to a constant multiple of
$H_{2a-l}(r) e^{-\frac{r^2}{2}}$ $(l=0,1)$ if $m=1$. 
\end{remNonumber}

\begin{remNonumber}
In \cite[Chapter 2]{xKo}, W. Myller-Lebedeff proved the following integral formula:
\begin{equation}%\label{eq:Myller-Lebedeff} \tag{4.2.5}
\label{eqn:myller}
  \int_{0}^{\infty} K_l(r,r';t)f_{a,l}(r')r'^{m-2}dr'=
  e^{-(a+\frac{m-1}2)t}f_{a,l}(r) \quad \text{for $a \ge l \ge 0$.}
\end{equation}
In view of \eqref{eqn:pilsum} and Proposition \ref{Lemma3.4} (2), 
the formula \eqref{4.2.1} in Theorem \ref{Theorem-C} implies
\eqref{eqn:myller} and vice versa.
The proof of \cite{xKo} is completely different from ours.
Here is a brief sketch:
For the partial differential operator
$$
  L:=\frac{\partial^2}{\partial x^2}+\frac{\alpha+1}x
  \frac{\partial}{\partial x}-\frac{1}x\frac{\partial}{\partial t}
  \quad \alpha >0,
$$
one has the following identity using Green's formula,
\begin{equation}\label{Green}
  \iint_D (v Lu- uL^\ast v)dtdx=\int_{\partial D}
(v \frac{\partial u}{\partial x}-u \frac{\partial v}{\partial x}+
\frac{\alpha+1}x uv)dt+\frac{1}x uv dx,
\end{equation}
for a domain $D\subset \mathbb R^2$. Here, $L^\ast$ denotes the
(formal) adjoint of $L$.
Now, we take $u(x,t):=t^n L_n^\alpha(\frac{x}t), v(x,t):=
(\frac{x}\xi)^\frac{\alpha}2
\frac{x}{\tau-t}e^{-\frac{x+\xi}{\tau-t}}
I_\alpha(\frac{2\sqrt{x\xi}}{\tau-t}), n\in \mathbb N,
\tau,\xi \in \mathbb R$
as solutions to $Lu=0, L^\ast v=0$ respectively, and
the domain $D$ as a rectangular domain
$D:=\set{(x,t)\in\mathbb{R}^2}{0<x<\infty, t_1<t<t_2}$ for some $t_1, t_2
\in \mathbb R$.
Then by the decay properties of $u$ and $v$,
the integrands in the right-hand side of \eqref{Green} 
vanish on $x=0$ and $x=\infty$.
Since the integral of the left-hand side of \eqref{Green}
vanishes,  the integral
$\int_{0}^{\infty}\frac{1}x u(x,t)v(x,t)dx$ becomes constant
with respect to $t$.
By taking the limit $t\to \tau$, we have
$\lim_{t \to \tau}\int_{0}^{\infty}\frac{1}x u(x,t)v(x,t)dx=
u(\xi,\tau)$ since $v(x,\tau)$ is proved to be a Dirac
delta function. Hence we obtain
$$
  \int_{0}^{\infty}\frac{1}x u(x,t)v(x,t)dx=u(\xi,\tau),
$$
which coincides with \eqref{eqn:myller} % (\ref{eq:Myller-Lebedeff})
by a suitable change of variables.
\end{remNonumber}

\subsection{Proof of Theorem \protect\ref{Theorem-C} (Case $\operatorname{Re}t=0$)}
\label{sec:4.3}
Suppose $t \in \sqrt{-1}\mathbb{R}$. 
Then, $\pi_{+,l}(e^{tZ})$ is a unitary operator on $\wL$. 
Suppose furthermore $t \notin 2\pi\sqrt{-1}\mathbb{Z}$. 
For $\varepsilon>0$ and 
$f \in \wL$, we have from Theorem \ref{Theorem-C} (1)
$$
  \pi_{+,l}(e^{(\varepsilon + t)Z})f=
  \int_0^\infty K_l^+(r, r'; \varepsilon+t)f(r'){r'}^{m-2}dr'. 
$$
By Proposition \ref{Lemma3.4} (1), the left-hand side converges to 
$\pi_{+,l}(e^{tZ})f$ in $\wL$ as $\varepsilon$ tends to $0$.
For the right-hand side, we have:

\begin{claim}
\label{clm:4.3.1}
For $t \in \sqrt{-1}\mathbb{R} \setminus 2\pi\sqrt{-1}\mathbb{Z}$,
\begin{alignat*}{2}
\lim_{\varepsilon \downarrow 0} 
\int_0^\infty K_l^+(r, r';\varepsilon+t) f_{a,l}(r')
{r'}^{m-2} dr'  
=  \int_0^\infty  K_l^+(r, r';t) f_{a,l}(r') {r'}^{m-2}dr'
\end{alignat*}
and the right-hand side converges absolutely.
\end{claim}

\begin{proof}
If $\varepsilon \in \mathbb{R}$ and
$t \in \sqrt{-1}\mathbb{R}$ then
$$
\Bigl|\sinh\frac{\varepsilon+t}{2}\Bigr|^2
= \Bigl|\sinh\frac{\varepsilon}{2}\Bigr|^2 + \Bigl|\sinh\frac{t}{2}\Bigr|^2
\ge \Bigl|\sinh\frac{t}{2}\Bigr|^2.
$$
Therefore, it follows from \eqref{eqn:Klupper} that
$$
|K_l^+ (r,r';\varepsilon+t)|
\le \frac{C(rr')^l}{|\sinh\frac{t}{2}|^{m-1+2l}}
$$
if $\varepsilon>0$ and 
$t \in \sqrt{-1}\mathbb{R} \setminus 2\pi\sqrt{-1}\mathbb{Z}$,
because $\varepsilon+t \in \Omega$ implies
$\alpha(\varepsilon+t) \ge 0$ and
$|\beta(\varepsilon+t)| < 1$.
Therefore, we have
$$
|K_l^+(r,r';\varepsilon+t) f_{a,l}(r') {r'}^{m-2}|
\le \frac{C r^l {r'}^{l+m-2}e^{-2r'}{|L_{a-l}^{m-2+4l}(4r')|}}
         {|\sinh\frac{t}{2}|^{m-1+2l}}.
$$
By the Lebesgue convergence theorem,
we have proved Claim.
\end{proof}

Since linear combinations of $f_{a,l}$ span a dense subspace of $\wL$, 2) is proved.  
\qed

\subsection{Weber's second exponential integral formula}
\label{sec:4.4}

{}From the semigroup law:
\begin{equation}
\label{4.7.1}
  \pi_+(e^{(t_1+t_2)Z})=\pi_+(e^{t_1Z})\circ \pi_+(e^{t_2Z})
  \quad  (\operatorname{Re}t_1, \operatorname{Re}t_2 >0),
\end{equation}
we get a representation theoretic proof of classical Weber's second exponential
integral for Bessel functions (see \cite[\S 13.31 (1)]{xWa}):

\begin{corSec}{\bf (Weber's second exponential integral)}
\label{Corollary-B}
Let $\nu$ be a positive integer, and $\rho, \alpha, \beta>0$.
We have the following integral formula 
\begin{equation}
\label{4.7.2}
\int_0^\infty e^{-\rho x^2} J_\nu(\alpha x)J_\nu(\beta x) x d x
= \frac{1}{2\rho} \exp\Bigl(-\frac{\alpha^2+\beta^2}{4\rho}\Bigr)
  I_\nu\Bigl(\frac{\alpha \beta}{2\rho}\Bigr).
\end{equation}
\end{corSec}
\begin{proof}

It follows from the semigroup law (\ref{4.7.1}) that the integral kernels for 
$\pi_+(e^{(t_1+t_2)Z})$ and $\pi_+(e^{t_1Z})\circ \pi_+(e^{t_2Z})$ must coincide.
Then, from Theorem \ref{Theorem-C}, we have 
\begin{equation}
\label{4.7.3}
  \int_0^\infty K_l^+(r, s; t_1)K_l^+(s, r'; t_2)s^{m-2}ds
  = K_l^+(r, r'; t_1+t_2).
\end{equation}
In view of (\ref{equation3}), the formula (\ref{4.7.2}) is obtained by 
(\ref{4.7.3}) by the change of variables,
\begin{alignat*}{2}
  &x= \sqrt{s}, 
  \quad \alpha= \frac{4e^\frac{\pi \sqrt{-1}}2 \sqrt{r}}{\sinh \frac{t_1}2 },
  \quad \beta= \frac{4e^\frac{\pi \sqrt{-1}}2 \sqrt{r'}}{\sinh \frac{t_2}2 },
  \quad \nu= m-2+2l, \\
  &\rho= 2\Bigl(\coth \frac{t_1}2 + \coth \frac{t_2}2  \Bigr).
\end{alignat*}
\end{proof}

\subsection{Dirac sequence operators} \label{Dirac sequence operators}

We shall state another corollary to Theorem \ref{Theorem-C}.
Let $\nu$ be a positive integer, $x, y \in \mathbb R, s \in \mathbb C$
such that $\operatorname{Re}s >0$. For a function $f$ on $\mathbb R$, 
let $\mathcal T_s$ be an operator defined by 
$$
  \mathcal T_s: f(x) \mapsto \int_0^\infty
  A(x, y; s) f(y)dy,
$$
with the kernel function 
\begin{equation} \label{definition of A}
  A(x,y ;s):= (xy)^\frac{1}2 
                    \frac{e^{-\frac{1}2(x^2+y^2) \coth s}}{\sinh s} I_\nu\Bigl(\frac{xy}{\sinh s}\Bigr).
\end{equation}
Then we have the following corollary. 
\begin{corSec}\label{Corollary-Dirac}
{\rm 1)}
The operators $\set{\mathcal T_s}{\operatorname{Re}s>0}$ form a semigroup of 
contraction operators on $L^2(\mathbb R_+, dx)$.

{\rm 2)} (Dirac sequence) $\lim_{s\to 0} || \mathcal T_s h- h ||_{L^2(\mathbb R_+, dx)}=0$
holds for all $h \in L^2(\mathbb R_+, dx)$. 
\end{corSec}

\begin{remNonumber}
For sufficiently small $s$, the semigroup $\set{\mathcal T_s}{\operatorname{Re}s>0}$
behaves like the Hermite semigroup (see \cite{xHo}) whose kernel is 
given by the following Gaussian (cf.\ \eqref{eqn:Mehler}):
$$
  \kappa (x, y; s)= 
  \frac{1}{\sqrt{2\pi \sinh s}} 
  e^{-\frac{x^2}2 \coth s +\frac{xy}{\sinh s}- \frac{y^2}2 \coth s}
$$
because $I_\nu (z) \sim \frac{1}{\sqrt{2\pi z}}e^z$ for sufficiently large
$z$ (see \cite[\S 7.23]{xWa} for the asymptotic behavior of $I_\nu(z)$). Note that 
it is stated in \cite[\S 5.5]{xHo} that the Hermite semigroup forms a \lq Dirac 
sequence'.
\end{remNonumber}

\begin{proof} Since we assume $\nu$ is a positive integer,
$\nu=m-2+2l$ for some $m >3$ and $l \in \mathbb Z$.
We change the variables
$r= \frac{x^2}4$, $r'= \frac{y^2}4$,  $t=2s$ and define a unitary map
\begin{equation} \label{Dirac-1}
\Phi: \wL \to L^2(\mathbb R_+, dx), \quad  (\Phi f)(x):=
\bigl(\frac{x}2\bigr)^\frac{2m-3}2 f\bigl( \frac{x^2}4 \bigr).
\end{equation}
Comparing (\ref{definition of A}) with (\ref{equation3}), we have
\begin{equation} \label{Dirac-2}
  \Phi^{-1} \circ \mathcal T_s \circ \Phi f= \pi_{+,l}(e^{2sZ})f
\end{equation}
Thus by Theorem \ref{Theorem-C} (1), 1) is proved.

2) We take a limit $s \downarrow 0$ of (\ref{Dirac-2}). By Proposition 
\ref{Lemma3.4} (1), the right-hand side equals $f$. Hence by putting $h:= \Phi f$,
we have $\lim_{s\downarrow 0} \mathcal T_s h=h$.
\end{proof}

\section{Integral formula for the semigroup}
\label{Main Theorem}

In this section, we shall give an explicit integral formula for the holomorphic semigroup
$\exp(t d \pi_+(Z)) = \pi_+(e^{t Z})$ on $\WL$
for $\operatorname{Re}t>0$, or more precisely, for $t \in \Omega=
\set{t \in \mathbb{C} }{ 
\operatorname{Re}t\geq 0 }\setminus 2\pi\sqrt{-1}\mathbb{Z}$
(see (\ref{definition of D})).
The main result of this section is Theorem \ref{Theorem-B}.
In particular, we give
 a proof of Theorem~A in Introduction.

\subsection{Result of the section}

Let $\langle x, x' \rangle$ be the standard inner product of 
$\mathbb R^m$, $|x|:=\sqrt{\langle x, x \rangle}$ be the norm.
We recall the notation from Subsection \ref{Differential operator D}:  
\begin{equation}
\label{eqn:psix}
\psi(x,x'):=2 \sqrt{2(|x||x'|+\langle x, x' \rangle)}
                =4|x|^\frac{1}2|x'|^\frac{1}2\cos \frac{\theta}2,
\end{equation} 
where $\theta\equiv \theta(x, x')$ is the angle between $x$ and $x'$ in $\mathbb R^m$.
Let us define a kernel function $\kakup{x}{x'}{t}$ 
on $\mathbb R^m \times \mathbb R^m \times \Omega$ 
by the following formula as in Introduction:
\begin{align}\label{5.2.1}
 K^+(x,x';t):= {}
&\frac{2^\frac{m-1}2 e^ {-2(|x|+|x'|) \coth \frac{t}{2}}   }
       {\pi^{\frac{m-1}{2}} \sinh^{\frac{m+1}{2}}\frac{t}{2}  }
 \psi(x, x')^{-\frac{m-3}{2}}
 I_{\frac{m-3}{2}  }
\Bigl( \frac{\psi(x, x')}{\sinh \frac{t}{2}} \Bigr)
\nonumber
\\
= {}
& \frac{2e^{-2(|x|+|x'|)\coth {\frac{t}{2}}}}
       {\pi^{\frac{m-1}{2}}\sinh^{m-1}\frac{t}{2}}
  \tilde{I}_{\frac{m-3}{2}}
  \Bigl( \frac{\psi(x,x')}{\sinh\frac{t}{2}} \Bigr),
\end{align}
where $I_\nu(z)$ is the modified Bessel function of the first kind 
and $\tilde{I}_\nu(z) := (\frac{z}{2})^{-\nu} I_\nu(z)$
is an entire function 
(see Subsection \ref{Bessel functions}). We note that $\sinh\frac{t}{2}$ in
the denominator is non-zero because $t \notin
2\pi\sqrt{-1}\mathbb{Z}$. 
Therefore,
$\kakup{x}{x'}{t}$ is a continuous function on 
$\mathbb{R}^m \times \mathbb{R}^m \times \Omega$. 

We recall from Proposition \ref{Lemma3.4} (3)
that $\pi_+(e^{tZ})$ is a contraction operator with operator norm 
$\| \pi_+(e^{tZ}) \|=e^{-\frac{m-1}2 \operatorname{Re}t}$.
Here is an integral formula of the holomorphic semigroup $\pi_+(e^{tZ})$:

\begin{thmSec}
[Integral formula for the semigroup]
\label{Theorem-B}

{\rm 1)} For $\operatorname{Re}t>0$, $\pi_+(e^{tZ})$ is a Hilbert--Schmidt 
operator on $\WL$, and is given by the following integral transform:
\begin{equation}
\label{5.2.2}  (\pi_+(e^{tZ})u)(x) = 
\int_{\mathbb R^m} \kakup{x}{x'}{t}u(x')\frac{dx}{|x|}
 \quad\text{for $u \in \WL$}.
\end{equation}
Here, the right-hand side converges absolutely.

{\rm 2)}  For  $t \in \sqrt{-1}\mathbb R$,
 $\pi_+(e^{tZ})$ is a unitary operator on $\WL$.
If $t \in \sqrt{-1}\mathbb{R}$ but 
$t \notin 2\pi\sqrt{-1}\mathbb{Z}$,
then the right-hand side of \eqref{5.2.2} converges absolutely for any 
$u \in L^1(\mathbb R^m, \frac{dx}{|x|})
\cap \WL$, in particular, for any $K$-finite vectors in 
$\WL$. Since $K$-finite vectors span a dense subspace of $\WL$, 
the integral formula \eqref{5.2.2} holds in the sense of $L^2$-convergence. 
\end{thmSec}

\begin{remNonumber}[Realization on the cone $C_+$]
Via the isomorphism
$L^2(C_+) \simeq L^2(\mathbb{R}^m, \frac{dx}{|x|})$
(see Subsection \ref{Schrodinger}), 
the above formula for $L^2(\mathbb{R}^m, \frac{dx}{|x|})$
can be readily transferred to the formula of the holomorphic extension
$\pi_+(e^{tZ})$ on $L^2(C_+)$.
For this, we define a continuous function 
$\widetilde{K^+}(\zeta, \zeta';t)$ 
on $C_+ \times C_+ \times \Omega$ 
by the following formula:
\begin{equation}
\label{eqn:Ckernel}
 \widetilde{K^+}(\zeta,\zeta'; t)
:= \frac{2e^ {-\sqrt{2}(|\zeta|+|\zeta'|) \coth \frac{t}{2}}   }
       {\pi^{\frac{m-1}{2}} \sinh^{m-1}\frac{t}{2}  }
 \tilde{I}_{\frac{m-3}{2}  }
\Bigl( \frac{2\sqrt{2\langle \zeta, \zeta' \rangle}}{\sinh \frac{t}{2}} \Bigr), 
\end{equation}
where $|\zeta|:=\sqrt{\langle \zeta, \zeta \rangle}=(\zeta_1^2+\cdots+\zeta_{m+1}^2)^\frac{1}2$.
Then, 
\begin{equation}
\label{eqn:Cint}  
(\pi_+(e^{tZ})u)(\zeta) = 
\int_{C_+} \widetilde{K^+}(\zeta, \zeta'; t)u(\zeta')d\mu(\zeta')
 \quad\text{for $u \in L^2(C_+)$}.
\end{equation}
\end{remNonumber}

\begin{remNonumber}[Weil representation]
For the Weil representation $\varpi$ of $G'$, the corresponding 
semigroup of contraction operators is the Hermite semigroup
$\set{\varpi(e^{tZ'})}{\operatorname{Re}t>0}$ (see Remark 3.6.3),
whose kernel function is given by
the Mehler kernel $\mathcal K(x,x'; t)$ (see Fact C in 
Subsection \ref{Our operator D and Hermite operator D};
see also \cite[\S5]{xHo}). 
\end{remNonumber}

The rest of this section is devoted to the proof of
 Theorem \ref{Theorem-B}.
Let us mention briefly a naive idea of the proof.
We observe that
 the action of 
$K(=SO(m+1) \times SO(2)\widetilde{}\:)$
on $\WL$ is hard to describe
because $K$ does not act on $\mathbb{R}^m$.
However, the action of
 its subgroup
$\R=SO(m)$ has a simple feature, that is,
 we have the following direct sum decomposition: 
\begin{equation}\label{5.2.3}
  \WL \simeq 
 \sideset{}{^\oplus}{\sum}_{l=0}^{\infty}
  \wL \otimes \Har{l}{m}.
\end{equation}
We have already proved in Theorem \ref{Theorem-C}
 that $K_l^+(r, r'; t)$ is the kernel of $\pi_{+,l}(e^{tZ})$ which is the restriction 
of $\pi_+(e^{tZ})$ in each $l$-component of the right-hand side of (\ref{5.2.3}).
Theorem \ref{Theorem-B} will be proved if we decompose $K^+$ into $K_l^+$.
This will be carried out in Lemma \ref{Lemma5.3}.
An expansion formula of $K^+$ by $K_l^+$ is not used in the proof of
 Theorem \ref{Theorem-B},
but might be of interest of its own.
We shall give it in Subsection \ref{subsec:5.7}.

\subsection{Upper estimates of the kernel function}
\label{Estimates of the kernel function}
In this subsection, we give an upper estimate
 of the kernel function $\kakup{x}{x'}{t}$. 
That parallels Lemma \ref{lem:Klupper}.
\begin{lemSec} \label{square-integrable}
Let $m \ge 2$.
{\rm 1)} There exists a constant $C > 0$ such that 
\begin{equation}\label{estimate of the kernel function}
  \bigl| \kakup{r\omega}{r'\omega'}{t}\bigr|\le 
  \frac{C}{ | \sinh \frac{t}2 |^{m-1} }
\, e^{-2\alpha(t) (1-|\beta(t)|) (r+r')},
\end{equation}
for any $r,r' \in \mathbb R_+$,
$\omega, \omega' \in S^{m-1}$, and 
$t \in \Omega$.
Here, 
$\alpha(t), \beta(t)$ are defined in \eqref{eqn:atxy} and
\eqref{eqn:btxy}, respectively.
  
{\rm 2)} If\/ $\operatorname{Re}t>0$, then 
$$
  \int_{\mathbb R^m} \int_{\mathbb R^m} \bigl| K^+(x, x'; t) \bigr|^2 
  \frac{dx}{|x|} \frac{dx'}{|x'|} < \infty.
$$

{\rm 3)} If\/ $\operatorname{Re}t>0$, then 
for a fixed $x \in \mathbb{R}^m$,
we have $K^+(x, \cdot\,; t) \in L^2(\mathbb{R}^m, \frac{dx'}{|x'|})$.

\end{lemSec}

\begin{proof}
By the definition \eqref{5.2.1} of
$K^+(x,x';t)$,
we have
$$
|K^+(r\omega,r'\omega';t)|
= \frac{2e^{-2(r+r')\coth\frac{t}{2}}}
       {\pi^{\frac{m-1}{2}}|\sinh\frac{t}{2}|^{m-1}}
\,\Bigl| \tilde{I}_{\frac{m-3}{2}}
    \left( \frac{\psi(r\omega,r'\omega')}{\sinh\frac{t}{2}} \right)
  \Bigr| .
$$
By \eqref{eqn:psix},
we have $|\psi(r\omega,r'\omega')| \le 4\sqrt{rr'}$.
Applying Lemma \ref{lem:estimate} with $\nu = \frac{m-3}{2}$,
we have
$$
|K^+(r\omega,r'\omega';t)|
\le \frac{2Ce^{-\alpha(t)(1-|\beta(t)|)(r+r')}}
         {\pi^{\frac{m-1}{2}} |\sinh\frac{t}{2}|^{m-1}}.
$$
Replacing $C$ with a new constant,
we get \eqref{estimate of the kernel function}.

The second and third statements follow from
 (\ref{estimate of the kernel function}) because
 $\alpha(t)>0$ and $|\beta(t)|<1$ if
$\operatorname{Re}t>0$.
\end{proof}

\subsection{Proof of Theorem \protect\ref{Theorem-B} (Case $\operatorname{Re}t>0$)}

Suppose $\operatorname{Re}t>0$.
We set
\begin{equation}\label{5.2.3.1}
  (S_tu)(x):=\int_{\mathbb R^m} K^+(x, x'; t) u(x') \frac{dx'}{|x'|}.
\end{equation}
By Lemma \ref{square-integrable} (2), 
we observe that $S_t$ is a Hilbert--Schmidt operator on $\WL$, and 
the right-hand side of (\ref{5.2.3.1}) converges absolutely for $u \in \WL$ by 
the Cauchy--Schwarz inequality and by 
Lemma \ref{estimate of the kernel function} (3).

The remaining assertion of Theorem \ref{Theorem-B} (1) is
 the equality $\pi_+(e^{tZ})=S_t$.
To see this, we observe from the definition \eqref{5.2.1} of
$K^+(x,x';t)$ that
$$
K^+(kx, kx';t)=K^+(x, x'; t) \quad\text{for all $k \in \R$}.
$$
Therefore, the 
operator $S_t$ intertwines the $\R$-action, and preserves each summand
of (\ref{4.1.1}). In light of the decomposition (\ref{eqn:pilsum}) 
of the operator $\pi_+(e^{tZ})$, the equality $\pi_+(e^{tZ})=S_t$
will follow from:

\begin{lemSec}\label{5.2.4}
Let $\operatorname{Re}t>0$. For every $l \in \mathbb N$, we have
\begin{equation}\label{5.2.4.1}
  \pi_{+,l}(e^{tZ}) \otimes \operatorname{id}=S_t|_{\wL \otimes \Har{l}{m}}.
\end{equation}
\end{lemSec}
We postpone the proof of Lemma \ref{5.2.4} until 
Subsection \ref{Proof of Lemma 5.2.4}.

Thus, the proof of Theorem \ref{Theorem-B} (1) is completed by
admitting Lemma \ref{5.2.4}.

\subsection{Proof of Theorem \protect\ref{Theorem-B} (Case $\operatorname{Re}t=0$)}

Suppose $\operatorname{Re}t=0$.
Then, by
Lemma \ref{square-integrable} (1), we have 
\begin{equation*}
  \bigl| K^+(x, x'; t)\bigr| \le  
  \frac{C}{ | \sinh \frac{t}2 |^{m-1} }
\end{equation*}
because $\alpha(t) = 0$. 
Therefore, the right-hand side of \eqref{5.2.2} converges absolutely
for any 
 $u \in L^1(\mathbb R^m, \frac{dx}{|x|}) \cap \WL$, 
as is seen by
\begin{equation*}
  \int_{\mathbb R^m} \bigl| K^+(x, x';t) u(x')\bigr| \frac{dx'}{|x'|} 
\le   \frac{C}{ | \sinh \frac{t}2 |^{m-1} }
     \int_{\mathbb R^m} | u(x') | \frac{dx'}{|x'|} <\infty.
\end{equation*}

By Proposition \ref{R_0 decomposition} (1),
we have
$\WL_K \subset L^1(\mathbb R^m, \frac{dx}{|x|})$. 
Hence, the right-hand side of \eqref{5.2.2} converges absolutely, 
in particular, for $K$-finite functions. 

Finally, let us show the last statement of (2).
Since $W_{a,l} (0\le l\le a)$ spans $\WL_K$ (see \eqref{R_0
decomposition}),
it is sufficient to prove 
\begin{equation}\label{eqn:pist}
\pi_+(e^{tZ}) = S_t
\quad\text{on ~$W_{a,l}~(0\le l \le a)$}.
\end{equation}
We recall
from Proposition \ref{branching-law} that every vector $u \in W_{a,l}$ 
is of the form
\begin{equation}
 u(r\omega)=f_{a,l}(r) \phi(\omega)
\end{equation}
for some $\phi \in \Har{l}{m}$(see \eqref{definition of f_a,l} for 
the definition).
Suppose $\varepsilon > 0$ and
$t \in \sqrt{-1} \mathbb{R} \setminus 2\pi \sqrt{-1} \mathbb{Z}$.
As in the proof of Claim \ref{clm:4.3.1}, 
we have
$$
|K^+ (x,x';\varepsilon+t)|
\le \frac{C'}{|\sinh\frac{t}{2}|^{m-1}}
$$
and therefore
\begin{equation*}
|K^+ (x,x';\varepsilon+t) f_{a,l}(r') \phi(\omega) {r'}^{m-2}|
\le
\frac{C {r'}^{m-2+l} e^{-2r'} |L_{a-l}^{m-2+2l} (4r')|}
     {|\sinh\frac{t}{2}|^{m-1}}.
\end{equation*}
Here,
$C := C' \max_{\omega\in S^{m-1}} |\phi(\omega)|$.
Hence, by the dominated convergence theorem,
we have
$$
\lim_{\varepsilon\downarrow 0} S_{\varepsilon+t} u 
= S_t u
$$
for any $u \in W_{a,l}$.

On the other hand,
by Proposition \ref{Lemma3.4} (1),
we have
$$
\lim_{\varepsilon\downarrow0} \pi_+(e^{(\varepsilon+t)Z}) u
= \pi_+(e^{tZ}) u.
$$
Since $\pi_+(e^{(\varepsilon+t)Z}) = S_{\varepsilon+t}$
for $\varepsilon > 0$,
we have now proved \eqref{eqn:pist}.

\subsection{Spectra of an $O(m)$-invariant operator}
\label{subsec:5.5}
The rest of this section is devoted to the proof of Lemma \ref{5.2.4}.

The orthogonal group $O(m)$ acts on $L^2(S^{m-1})$ as a unitary
representation, 
and decomposes it into irreducible representations as follows:
$$
L^2(S^{m-1}) \simeq \sideset{}{^\oplus}\sum_{l=0}^\infty
\mathcal{H}^l(\mathbb{R}^m).
$$
Since this is a multiplicity-free decomposition,
any $O(m)$-invariant operator $S$ on $L^2(S^{m-1})$ acts on each
irreducible component as a scalar multiplication. 

The next lemma gives an explicit formula of the spectrum 
for an $O(m)$-invariant integral operator on $C(S^{m-1})$ in the
general setting. 
This should be known to experts, but for the 
convenience of the readers, we present it in the following form:
\begin{lemSec}
\label{Lemma5.4}
For a continuous function $h$ on the closed interval $[-1, 1]$,
we consider the following integral transform:
\begin{equation}
\label{definition of S_h}
  S_h :L^2(S^{m-1}) \to L^2(S^{m-1}),~~\phi(\omega)  \mapsto 
   \int_{S^{m-1}} h(\langle \omega, \omega' \rangle)
   \phi(\omega') d\omega' .
\end{equation}
Then, $S_h$ acts on $\mathcal{H}^l(\mathbb{R}^m)$ by a 
scalar multiplication of $c_{l,m}(h)\in \mathbb{C}.$
The constant $c_{l,m}(h)$ is given by  
\begin{equation}
\label{5.4.1}
          c_{l,m}(h) = \frac{2^{m-2}\pi^{\frac{m-2}{2}}l!}
                            {\Gamma(m-2+l)}
          \int_0^\pi h(\cos\theta) \widetilde{C}_l^{\frac{m-2}{2} }(\cos\theta)
          \sin^{m-2}\theta d\theta,
\end{equation}
where $\widetilde{C}_l^{\frac{m-2}{2}  }(x)$ denotes the normalized Gegenbauer polynomial 
(see \eqref{eqn:nGegen}). 
\end{lemSec}
\begin{example}[see {\cite[Introduction, Lemma 3.6]{xHe}}]
For $h(x):=e^{\sqrt{-1}\lambda x}$, $c_{l,m}(h)$ amounts to
$$
  c_{l,m}(h)=(2\pi)^\frac{m}2 e^{\frac{\sqrt{-1}}{2}\pi l } 
 \lambda^{-\frac{m-2}2}J_{\frac{m-2}2+l}(\lambda).
$$
\end{example}

\begin{example}
\label{ex:coefIB}
We set 
$\tilde{I}_\nu(z) = \left(\frac{z}{2}\right)^{-\nu} I_\nu(z)$
(see \eqref{eqn:Itilde}).
For $h(s):= \tilde{I}_{\frac{m-3}{2}} (\alpha\sqrt{1+s})$,
we have
\begin{equation}
\label{eqn:coefIB}
c_{l,m}(h) = 2^{\frac{3m-4}{2}} \pi^{\frac{m-1}{2}}
  \alpha^{-m+2} I_{m-2+2l} (\sqrt{2} \alpha)
\end{equation}
\end{example}

\begin{proof}[Proof of Example \ref{ex:coefIB}]
We apply Lemma \ref{Lemma5.5} with $\nu = \frac{m-3}{2}$. 
Then, we have
\begin{align*}
&\int_0^\pi \tilde{I}_{\frac{m-3}{2}} (\alpha\sqrt{1+\cos\theta})
 \widetilde{C}_l^{\frac{m-2}{2}} (\cos\theta) \sin^{m-2} \theta \, d\theta
\\
&= \frac{2^{\frac{m}{2}} \sqrt{\pi}\, \Gamma(m-2+l)}
        {\alpha^{m-2} l!}
 \, I_{m-2+2l} (\sqrt{2} \alpha).
\end{align*}
Hence, \eqref{eqn:coefIB} follows from Lemma \ref{Lemma5.4}.
\end{proof}

\begin{proof}[Proof of Lemma \ref{Lemma5.4}]
1) The operator $S_h$ intertwines the $O(m)$-action because
$h(\langle k\omega, k\omega' \rangle)=h(\langle \omega, \omega' \rangle)$ 
for $k \in O(m)$. Hence it follows from Schur's lemma that 
$S_h$ acts on each irreducible $O(m)$-subspace
 $\mathcal{H}^l(\mathbb{R}^{m})$ by the multiplication of a constant,
which we shall denote by $c_{l, m}(h)$ for $l=0,1,2,\ldots$.
Thus, we have
\begin{equation} \label{5.4.2}
  (S_h \phi)(\omega)=c_{l,m}(h) \phi(\omega)
\quad\text{for $\phi \in \mathcal{H}^l(\mathbb{R}^m)$}.
\end{equation}

To compute the constant $c_{l,m}(h)$,  we use the following coordinate:
$$
[0,\pi )\times S^{m-2} \to S^{m-1}, \quad 
(\theta, \eta) \mapsto \omega =(\cos\theta, \sin\theta\cdot \eta).
$$ 
With this coordinate, we have $d\omega=\sin^{m-2}\theta d\theta d\eta$. 

We set $\omega_0=(1,0, \cdots, 0)$. 
Now, we take
$\phi(\omega):= \widetilde{C}_l^{\frac{m-2}{2}}(\langle\omega,\omega_0\rangle) 
\in \mathcal{H}^l(\mathbb{R}^{m})$, 
which is an $O(m-1)$-invariant spherical harmonics. 
Then the equation (\ref{5.4.2}) for $\omega=\omega_0$ 
amounts to 
\begin{equation*} 
\frac{2\pi^{\frac{m-1}{2}}}{\Gamma(\frac{m-1}2)}
   \int_0^\pi h(\cos\theta) \widetilde{C}_l^{\frac{m-2}{2}}(\cos\theta)
   \sin^{m-2}\theta d\theta 
= \frac{\sqrt{\pi}\,\Gamma(m-2+l)}{2^{m-3}l!\, \Gamma(\frac{m-1}{2})} c_{l,m}(h),
\end{equation*}
because $\operatorname{vol}(S^{m-2})
=\frac{2\pi^{\frac{m-1}{2}}}{\Gamma(\frac{m-1}2)}$ and
$\phi(\omega_0)=\widetilde{C}_l^{\frac{m-2}{2}}(1)
=\frac{\sqrt{\pi}\,\Gamma(m-2+l)}{2^{m-3}l!\, \Gamma(\frac{m-1}{2})}$
(see \eqref{eqn:Gegen1}). 
Hence, (\ref{5.4.1}) is proved.
\end{proof}

\subsection{Proof of Lemma \protect\ref{5.2.4}}
\label{Proof of Lemma 5.2.4}
This subsection gives a proof of Lemma \ref{5.2.4}.

We recall from Theorem \ref{Theorem-C} (1) that the kernel function of $\pi_{+,l}(e^{tZ})$
is given by $K_l^+(r, r';t)$ (see (\ref{equation3}) for definition). 
Therefore, the equation (\ref{5.2.4.1})
is equivalent to the following equation between kernel functions: 
\begin{lemSec}
\label{Lemma5.3}
For $\phi \in \Har{l}{m}$, we have \begin{equation}
\label{5.3.2}
  \frac{1}2 \int_{S^{m-1}} \kakup{r\omega}{r'\omega'}{t} 
  \phi(\omega')d\omega'
  = K_l^+(r, r'; t)\phi(\omega).
\end{equation}
\end{lemSec}

\begin{proof}[Proof of Lemma \ref{Lemma5.3}]
We set
\begin{equation}
\label{eqn:hrrts}
h(r,r',t,s) := \frac{2e^{-2(r+r')\coth\frac{t}{2}}}
                    {\pi^{\frac{m-1}{2}} \sinh^{m-1} \frac{t}{2}}
               \, \tilde{I}_{\frac{m-3}{2}}
               \left( \frac{2\sqrt{2rr'(1+s)}}{\sinh\frac{t}{2}}
                    \right).
\end{equation}
By the definition \eqref{5.2.1} of $K^+(x,x';t)$,
we have
\begin{equation*}
K^+(r\omega,r'\omega';t) = h(r,r',t,\langle\omega,\omega'\rangle).
\end{equation*}
Consider now the integral transform 
$S_{h(r,r',t,\cdot)}$ on $L^2(S^{n-1})$ with kernel
$h(r,r',t,\langle\omega,\omega'\rangle)$.
Then it follows from Example \ref{ex:coefIB} with
$\alpha = \frac{2\sqrt{2rr'}}{\sinh\frac{t}{2}}$
that
\begin{equation}
\label{eqn:KtoKl}
c_{l,m}(h(r,r',t,\cdot)) = 2K_l^+(r,r';t).
\end{equation}
Then, \eqref{5.3.2} is a direct consequence of Lemma \ref{Lemma5.4}. 
\end{proof}

\subsection{Expansion formulas}
\label{subsec:5.7}

We recall that the kernel functions for the semigroups $\pi_+(e^{tZ})$ and 
$\pi_{+,l}(e^{tZ})$ are given by $K^+(r\omega, r'\omega'; t)$ and 
$K_l^+(r, r'; t)$. In this subsection, we shall give expansion formulas for
$K^+(r\omega, r'\omega'; t)$ arising from the decomposition 
(see \eqref{eqn:pilsum})
\begin{equation*} 
  \pi_+(e^{tZ})= \sideset{}{^\oplus}{\sum}_{l=0}^\infty
  \pi_{+,l}(e^{tZ})\otimes \operatorname{id}.
\end{equation*}
\begin{proposition}[Expansion formulas] \label{prop:expansion}
Let $m>1$. 

1) The kernel function $K^+(x,x';t)$ (see \eqref{5.2.1}) 
has the following expansion:
\begin{equation}\label{eqn:expansion1}
  K^+(r\omega, r'\omega';t)= 
     \frac{1}{\pi^\frac{m}2}\sum_{l=0}^\infty (\frac{m-2}2+l)
      K_l^+(r, r';t) \widetilde{C}_l^\frac{m-2}2(\langle \omega, \omega' \rangle).
\end{equation}

2) The special value $t=\pi \sqrt{-1}$ for \eqref{eqn:expansion1} 
yields the expansion formula for the Bessel function:
\begin{equation}\label{eqn:expansion1-1}
 \tilde{J}_{\nu-\frac{1}2}(\sqrt z \cos \frac{\theta}2)=
 \frac{2^{4\nu}}{\sqrt \pi}
 \sum_{l=0}^\infty (\nu+l) (-1)^l 
 \frac{J_{2\nu+2l}(\sqrt z)}{z^\nu} \widetilde{C}_l^\nu(\cos \theta)
\end{equation}
for $z \in \mathbb R_+, \nu \in \frac{1}2 \mathbb Z$.
\end{proposition}
\begin{remNonumber}[Weil representation, Gegenbauer's expansion]
\label{rem:5.7.2}
Let us compare the above result with
 the case of the Weil representation of $G'$.
Then, by a similar argument to the proof of Proposition
 \ref{prop:expansion},
we can show that
the Mehler kernel $\mathcal K$ 
(see Subsection \ref{Our operator D and Hermite operator D})
has the following decomposition: 
\begin{equation}
\label{eqn:expansion2}
  \mathcal K(r\omega, r'\omega';t)=
      \frac{1}{2\pi^\frac{m}2}\sum_{l=0}^\infty (\frac{m-2}2+l)
      \mathcal K_l(r, r';t) \widetilde{C}_l^\frac{m-2}2(\langle \omega, \omega' \rangle)
\end{equation}
if $m>1$. 

In light of the formula (see Subsection
\ref{Our operator D and Hermite operator D}) % 1.2)
$$
  \mathcal K(r\omega, r'\omega'; \pi \sqrt{-1})=
  \frac{1}{(2\pi \sqrt{-1})^\frac{m}2} e^{-\sqrt{-1}rr'\langle \omega, \omega' \rangle},
$$
the special value at $t=\pi \sqrt{-1}$ for \eqref{eqn:expansion2} yields the following
expansion formula for the exponential function known as 
{\it Gegenbauer's expansion}
(\cite{xGegen}, see also \cite[Chapter~XI, \S 11.5]{xWa}):
\begin{equation}
  e^{\sqrt{-1}z \cos \phi}= 2^\nu 
  \sum_{m=0}^\infty (\nu+m) \sqrt{-1}^m \frac{J_{\nu+m}(z)}{z^\nu}
  \widetilde{C}_m^\nu(\cos \phi)
\end{equation}
for $z \in \mathbb R, \nu \in \frac{1}2 \mathbb Z$. 
This formula corresponds to \eqref{eqn:expansion1-1}.

If $m=1$, we have
\begin{equation}\label{eqn:expansion3}
  \mathcal K(x, x'; t)= \frac{1}2 \bigl(\mathcal K_0 (|x|, |x'|; t) \operatorname{id}(x) +
\mathcal K_1(|x|, |x'|; t) \operatorname{sgn}(x) \bigr),
\end{equation}
which is immediately verified 
by the formulas for $\mathcal K_0, \mathcal K_1$ and 
$\mathcal K$ (see Remark 4.1.2).
\end{remNonumber}

\begin{proof}[Proof of Proposition \ref{prop:expansion}]
1) First we prepare a lemma.

\begin{lemSec}\label{lem:expansion}
Let $f \in C[-1,1]$. Then $f$ is expanded into the infinite series:
\begin{equation}
  f(x)=\frac{1}{2\pi^\frac{m}2}
          \sum_{l=0}^\infty \bigl( \frac{m-2}2+l \bigr) c_{l,m}(f) 
        \widetilde{C}_l^\frac{m-2}2(x),
\end{equation}
where $c_{l,m}(f)$ is the constant defined in \eqref{5.4.1}.
\end{lemSec} 
\begin{proof}[Proof of Lemma \ref{lem:expansion}]
Applying the expansion formula
$f(x) = \sum_{l=0}^\infty \alpha_l^\nu (f) \widetilde{C}_l^\nu(x)$
(see \eqref{eqn:Gegenexpansion}) with
$\nu = \frac{m-2}{2}$,
 we have
\begin{equation}\label{eqn:lemma2}
  \alpha_l^{\frac{m-2}{2}}(f)
  = \frac{(\frac{m-2}2+l)}{2\pi^\frac{m}2}
  c_{l,m}(f),
\end{equation}
from \eqref{5.4.1} and \eqref{eqn:coeffGegen}. 
Hence, we get the lemma.
\end{proof}

Let $h(r,r',t,s)$ be as in \eqref{eqn:hrrts}.
We recall from \eqref{eqn:KtoKl}
$$
   c_{l,m}(h(r,r',t,\cdot))=2K_l^+(r, r';t).
$$
Therefore, by Lemma \ref{lem:expansion}, we get (\ref{eqn:expansion1}).

2) Substitute $t=\pi \sqrt{-1}$ and put
$$
 z:=16rr', \quad \cos \theta:=\langle \omega, \omega' \rangle, \quad 
 \nu=\frac{m-2}2,
$$
we get (\ref{eqn:expansion1-1}).
\end{proof}

\section{The unitary inversion operator}
\label{Integral formula for the inversion operator.}

\subsection{Result of the section}
\label{sec:6.1}
We define the \lq inversion' element $w_0 \in G$ by 
$$
 w_0:=e^{\pi\sqrt{-1}Z}.
$$

Then, clearly, $w_0$ has the following properties:

1) $w_0$ is of order four.

2) $w_0$ normalizes $\Mmax A$ and
$
\operatorname{Ad}(w_0)\mathfrak{n}^\text{max}=\overline{\mathfrak{n}^\text{max}}. 
$

3) The group $G$ is generated by $\overline{\Pmax}$ and $w_0$.

We note that if $m$ is odd then $e^{\pi\sqrt{-1}Z}$ is equal to 
$\begin{pmatrix}  I_{m+1}  & 0  \\ 0  &  -I_2 \end{pmatrix}$
in $SO_0(m+1,2) = G/\{1,\eta\}$.

This section gives an explicit integral formula
of the unitary operator $\pi_+(w_0)$ on $\WL$.
In light of $w_0 = e^{\pi\sqrt{-1}Z}$ and $\pi \sqrt{-1} \in \Omega$, 
we define the following kernel functions by substituting $t=\pi\sqrt{-1}$ into 
(\ref{5.2.1}) and (\ref{equation3}), respectively: 
\begin{alignat}{2}
\label{def_of_Kplus}
K^+(x, x'):={} & \kakup{x}{x'}{\pi\sqrt{-1}}
\nonumber
\\
={}& 
               \frac{2}{e^{\frac{m-1}{2}\pi \sqrt{-1}}\pi^{\frac{m-1}{2}} }
               \sqrt{2\langle \zeta, \zeta' \rangle }^{-\frac{m-3}{2}}
               J_{\frac{m-3}{2}}(2\sqrt{2\langle \zeta, \zeta' \rangle}), \\
   K_l^+(r, r'):={} & K_l^+(r,r';\pi\sqrt{-1})
\nonumber
\\
={}&
               2(-1)^l e^{-\frac{m-1}{2}\pi \sqrt{-1}} (rr')^{-\frac{m-2}{2}}
               J_{m-2+2l}(4\sqrt{rr'}).
\end{alignat}  
Then, the following result is a direct consequence of Theorems
\ref{Theorem-B} and \ref{Theorem-C}. 

\begin{thmSec}[Integral formula for the unitary inversion operator]
\label{thm:D}
The unitary operator $\pi_+(w_0): \WL \to \WL$
is given by the following integral transform:
\begin{equation}
\label{6.1.1}
(Tu)(x):=\int_{\mathbb R^m} K^+(x, x')u(x')\frac{dx}{|x|}, \quad u \in \WL.
\end{equation}
Here, the integral \eqref{6.1.1}  
converges absolutely for any $u \in L^1(\mathbb R^m, \frac{dx}{|x|}) \cap \WL$, 
in particular, for any $K$-finite vector, 
hence the equation \eqref{6.1.1} holds in the sense of $L^2$-convergence.
\end{thmSec}

The substitution of $t = \pi\sqrt{-1}$ into \eqref{eqn:pilsum} gives
the decomposition
$$
\pi_+(w_0) = \sideset{}{^\oplus}{\sum}_{l=0}^\infty \pi_{+,l}
             (w_0) \otimes \operatorname{id}.
$$
As $\pi_+(w_0)$ is given by the kernel $K^+(x,x')$,
so is $\pi_{+,l}(w_0)$ by $K_l^+(r,r')$.
Thus, we have:

\begin{thmSec}[Radial part of the unitary inversion operator]\label{thm:E}
The unitary operator
$\pi_{+,l}(w_0):
 L^2(\mathbb{R}_+, r^{m-2} dr) \to L^2(\mathbb{R}_+, r^{m-2} dr)$
is given by the following integral transform:
\begin{equation}
\label{equation5}
(T_lf)(r):=\int_0^\infty K_l^+(r, r')f(r'){r'}^{m-2}dr'.
\end{equation}  
Here, the integral \eqref{equation5}  
converges absolutely for any $f \in L^1(\mathbb R_+, r^{m-2}dr) \cap \wL$.
 In particular, 
 the equation \eqref{equation5} holds in the sense of $L^2$-convergence.

\end{thmSec}

\begin{remNonumber}[Weil representation]
In the case of the Weil representation $\varpi$ of $G'$, 
the counterparts of Theorem \ref{thm:D} and Theorem \ref{thm:E} 
can be stated as follows:
let $\omega_0 := e^{\pi\sqrt{-1}Z'}$.

We define kernel functions $\mathcal K$ and $\mathcal K_l$ by the formulas
\begin{alignat*}{1}
  \mathcal K(x, x'):=&\mathcal K(x, x'; \pi \sqrt{-1})=
\frac{1}{(2\pi\sqrt{-1})^\frac{m}2} e^{-\sqrt{-1} \langle x, x' \rangle},  \\
  \mathcal K_l(r, r'):=&\mathcal K_l(r, r'; \pi \sqrt{-1})=
  \sqrt{-1}^{-\frac{m-1+2l}2}(rr')^{-\frac{m-2}2} 
  J_\frac{m-2+2l}2 (rr'). 
\end{alignat*}
Then, 

1) The unitary operator $\varpi(w_0): L^2(\mathbb R^m) \to 
L^2(\mathbb R^m)$ is given by 
\begin{equation*}
  (\varpi(w_0)u)(x)=\int_{\mathbb R^m} \mathcal K(x,x')u(x')dx'.
\end{equation*}
Hence we see that $\varpi(w_0)$ is nothing but the Fourier transform.

2) The \lq radial' part of $\varpi(w_0)$, namely, 
the unitary operator $\varpi_l(w_0):L^2(\mathbb R_+, r^{m-1}dr) \to 
L^2(\mathbb R_+, r^{m-1}dr)$ (see Remark 4.1.2) is given by
\begin{equation*}
 (\varpi_l(w_0)f)(r)=\int_0^\infty \mathcal K_l(r, r')f(r'){r'}^{m-1}dr'.
\end{equation*}
\end{remNonumber}

\subsection{Inversion and Plancherel formula}

It follows from \eqref{eigenfunction_of_T} that
$\pi(w_0)^2 = \pi_+(e^{2\pi\sqrt{-1}Z}) = (-1)^{m+1}\operatorname{id}$. 
Thus,
as an immediate consequence of Theorem \ref{thm:D}, we have:

\begin{corSec}
\label{cor:E}
The integral transform 
\begin{alignat*}{1}
T&: u(x) \mapsto \int_{\mathbb R^m} K^+(x,x') \frac{u(x')}{|x'|}dx', \\
  K^+(x, x')&:=\frac{2^\frac{m-1}2 \psi(x,x')^{-\frac{m-3}2}}
           {e^{\frac{m-1}2 \pi \sqrt{-1}} \pi^\frac{m-1}2 } 
  J_\frac{m-3}2 (\psi(x, x'))
\end{alignat*} 
is a unitary operator on $L^2(\mathbb R^m, \frac{dx}{|x|})$
of order two ($m$: odd) and of order four ($m$: even), that is, we have:
\begin{alignat*}{1}
&\text{\rm (Inversion formula)} \quad T^{-1} = (-1)^{m+1} T,\\
&\text{\rm (Plancherel formula)} \quad 
\| Tu \|_{L^2(\mathbb R^m, \frac{dx}{|x|})} =\|u\|_{L^2(\mathbb R^m, \frac{dx}{|x|})}
\quad \text{for all}~ u \in \WL.
\end{alignat*}
\end{corSec}

\begin{remNonumber}
[{Weil representation, see \cite[Corollaries 5.7.3 and 5.7.4]{xHo})}]
\label{rem:Weil Fourier}
Let us compare Corollary \ref{cor:E} with the corresponding result for
the Schr\"{o}dinger model $L^2(\mathbb{R}^m)$ of the Weil
representation $\varpi$.
The unitary operator $\varpi(w_0)$ corresponding to the ``inversion
element'' $w_0$ is given by the (ordinary) Fourier transform
$\mathcal{F}$ (see Fact \ref{fact:D} 
in Subsection \ref{Our operator D and Hermite operator D}).
As is well-known,
$\mathcal{F}$ is a unitary operator of order four.
This reflects the fact that $w_0^4 = e$ in $\Mp$.
\end{remNonumber}

\subsection{The Hankel transform}
Similar to Corollary \ref{cor:E}, from the equality 
$\pi_+(w_0)^2= (-1)^{m+1} \text{id}$ 
and Theorem \ref{thm:E}, we have:
\begin{corSec}\label{cor:F}
Let $\nu$ be a positive integer, and $x, y \in \mathbb R$. Then the integral transform
$$
  \mathcal T_\nu: f(x) \mapsto \int_0^\infty J_\nu(xy)f(y) \sqrt{xy} dy
$$
is a unitary operator on $L^2((0, \infty), dx)$ of order two. Hence we have:
\begin{alignat*}{1}
&\text{\rm (Inversion formula)} ~~~  \mathcal T_\nu^{-1}
=\mathcal T_\nu,
\\ 
 &\text{\rm (Plancherel formula)}  ~~~
\| \mathcal T_\nu f \|_{L^2(\mathbb R_+, dx)} = \|f\|_{L^2(\mathbb R_+, dx)}.
\end{alignat*}
\end{corSec}
\begin{remNonumber} The unitary operator $\mathcal T_\nu$ coincides with the 
{\it Hankel transform} (see \cite[Chapter VIII]{xerdInt}) and the property
$\mathcal T_\nu^2=\operatorname{id}$ in the corollary 
corresponds to its classically known {\it reciprocal formula} 
due to Hankel \cite{xhankel} (see also
\cite[\S 8.1 (1)]{xerdInt}, \cite[\S 14.3 (3)]{xWa}).
The Parseval-Plancherel formula for the Hankel transform 
goes back to Macaulay-Owen \cite{xMa}.
\end{remNonumber}

\begin{proof}[Proof of Corollary \ref{cor:F}]
Since we assume $\nu$ is a positive integer,
$\nu=m-2+2l$ for some $m >3$ and $l \in \mathbb Z$.
We change the variables
$r= \frac{x^2}4$, $r'= \frac{y^2}4$ and define a unitary map
$\Phi: \wL \to L^2(\mathbb R), (\Phi f)(x):=
\bigl(\frac{x}2\bigr)^\frac{2m-3}2 f\bigl( \frac{x^2}4 \bigr)$.
Then, we have
$\mathcal{T}_\nu = e^{-(l+\frac{m-1}{2})\pi\sqrt{-1}}\Phi\circ T_l
 \circ \Phi^{-1}$.
Since $T_l = \pi_{+,l}(w_0)$ by Theorem \ref{thm:E},
$T_l$ has its inverse as
$T_l^{-1} = (-1)^{m+1} T_l$,
$\mathcal{T}_\nu$ has its inverse given by
$$
\mathcal{T}_\nu^{-1} = e^{(l+\frac{m-1}{2})\pi\sqrt{-1}}
 \Phi \circ (-1)^{m+1} T_l \circ \Phi^{-1}
= \mathcal{T}_\nu.
$$
Hence,
 the corollary follows.
\end{proof}

\subsection{Forward and backward light cones}
So far, we have discussed only the irreducible unitary representation 
$\pi_+$ realized on $L^2(C_+)$ for the forward light cone.
In this subsection, let us briefly comment on
$(\pi_-, L^2(C_-))$ for the backward light cone
and $(\pi, L^2(C))$ for $C=C_+ \cup C_-$. 

Since $d\pi_-(Z)=-d\pi_+(Z)$
in the polar coordinate representation, we can define a semigroup
of contraction operators $\set{\pi_-(e^{tZ})}{\operatorname{Re}t <0}$ 
similarly on $L^2(C_-)$.
All the statements about the semigroup $\pi_+(e^{tZ})$ also hold by changing 
the signature $t \to -t$ and replacing $C_+$ by $C_-$.
  
We define a function $\kaku$ on $C \times C$ by
\begin{equation}\label{eqn:KC}
   \kaku := 
   \begin{cases}
     K^+(\zeta, \zeta')   & \text{if} ~~ \langle \zeta, \zeta' \rangle \ge 0   \\
     0      & \text{if}~~ \langle \zeta, \zeta' \rangle <0.
   \end{cases}
\end{equation}
Here we note:
$$\set{(\zeta, \zeta') \in C \times C}{\langle \zeta, \zeta' \rangle \ge 0}
= (C_+ \times C_+) \cup (C_- \times C_-).
$$
Originally, $K^+(\zeta, \zeta')$ was defined on $C_+ \times C_+$ 
(see (\ref{def_of_Kplus})). Since $K^+(\zeta, \zeta')$ depends only on 
the inner product $\langle \zeta, \zeta' \rangle$, we can define 
$K^+(\zeta, \zeta')$ also on $C_- \times C_-$.
\begin{corSec}
\label{cor:6.4.1}
The unitary operator $\pi(w_0): L^2(C) \to L^2(C)$
coincides with the integral transform defined by
\begin{equation}
\label{thm3}
T: L^2(C) \to L^2(C),\quad 
u \mapsto \int_{C} K(\zeta, \zeta')u(\zeta')d\mu(\zeta').
\end{equation}
\end{corSec}

\begin{remNonumber}
We note that the kernel function $\kaku$ is supported on the proper subset
$\set{(\zeta, \zeta') \in C \times C}{\langle \zeta, \zeta' \rangle \ge 0}$ 
of $C \times C$. More generally, in the case of the minimal representation of $O(p,q)$
with $p+q: \text{even}, \ge 8$, we shall see in \cite{xkmano4} 
that the integral kernel $K_{p,q}(\zeta, \zeta')$ representing 
the inversion element is also supported on the proper subset
$\set{(\zeta, \zeta') \in C \times C}{\langle \zeta, \zeta' \rangle \ge 0}$ 
if both integers $p, q$ are even. This feature fails if both $p, q$ are odd.
In fact, the support of
$K_{p,q}$ is the whole space $C \times C$ then.
\end{remNonumber}

\section{Explicit actions of the whole group on $L^2(C)$}
\label{sec:formulaG}

Building on the explicit formula of 
$\pi(w_0) = \pi(e^{\pi\sqrt{-1}Z})$ on
$L^2(C)$
(see Corollary \ref{cor:6.4.1})
and that of $\pi(g)$ $(g \in \overline{\Pmax})$
(see Subsection \ref{Schrodinger model of the minimal representation.}),
we can find an explicit formula of the minimal representation for the
whole group $G$.

For simplicity, this section treats the case where $m$ is odd.
The main result of this section is Theorem \ref{thm:formulaG} for the
action of $O(m+1,2)$.

\subsection{Bruhat decomposition of $O(m+1,2)$}
\label{sec:BruhatG}

We recall the notation in Subsection \ref{Notations and definitions.}.
In particular,
 $\overline{\Pmax} = \Mmax A \overline{\Nmax}$ 
is a maximal parabolic subgroup of
$O(m+1,2)$.
Hence, $O(m+1,2)$ is expressed as the disjoint union:
\begin{align*}
O(m+1,2)
&= \overline{\Pmax} \amalg \overline{\Pmax} w_0 \overline{\Pmax}
\\
&= \overline{\Pmax} \amalg \overline{\Nmax} A\Mmax w_0
   \overline{\Nmax}. 
\end{align*}
We begin by finding $a,b\in\mathbb{R}^{m+1}$, $t\in\mathbb{R}$ and
$m\in\Mmax$ such that
\begin{equation}
\label{eqn:gBruhat}
g = \nbar{b} e^{tE} m w_0 \nbar{a}
\in \overline{\Nmax} A \Mmax w_0 \overline{\Nmax}
\end{equation}
holds for $g \notin \overline{\Pmax}$.

Suppose $g$ is of the form \eqref{eqn:gBruhat},
and we write $m=\delta m_+$ $(\delta=\pm1, m_+\in\Mmax_+)$.
Then, in light of the formulas \eqref{eqn:nbn}--\eqref{eqn:a+m},
we have
\begin{align}
\label{eqn:geQ}
g(e_0 - e_{m+2})
&= \nbar{b} e^{tE}(\delta m_+)w_0\nbar{a}(e_0 - e_{m+2})
\nonumber
\\
&= \nbar{b} e^{tE}(\delta m_+)w_0(e_0 - e_{m+2})
\nonumber
\\
&= \nbar{b} e^{tE}(\delta m_+)w_0(e_0 - e_{m+2})
\nonumber
\\
&=\delta e^t\nbar{b}(e_0 + e_{m+2})
\nonumber
\\
&=\delta e^t\, \trans (1-Q(b), 2b_1,\dots,2b_{m+1}, 1+Q(b) ).
\end{align}

Now, we set
\begin{equation}
\label{eqn:xge}
x \equiv \trans (x_0,\dots,x_{m+2}) := g(e_0 - e_{m+2}).
\end{equation}
Then, solving \eqref{eqn:geQ}, we have
\begin{align*}
& x_0 + x_{m+2} = 2 \delta e^t,
\\
& \frac{x_j}{x_0+x_{m+2}} = b_j
  \quad (1 \le j \le m+1).
\end{align*}
Likewise, if we set
\begin{equation}
\label{eqn:yge}
y \equiv \trans (y_0,\dots,y_{m+2})
:= g^{-1} (e_0 - e_{m+2}),
\end{equation}
the following relation:
$$
\frac{y_j}{y_0+y_{m+2}} = -a_j
\quad (1 \le j \le m+1)
$$
must hold, because
$g^{-1} = \nbar{-a} e^{tE} (w_0 m^{-1} w_0^{-1}) w_0
 \nbar{-b}$.

Let $\langle\ ,\ \rangle$ denote the standard positive definite inner
product on $\mathbb{R}^{m+3}$. 
It follows from $g \in O(m+1,2)$ that
\begin{align*}
x_0 + x_{m+2}
&= \langle e_0 + e_{m+2}, g(e_0 - e_{m+2}) \rangle
\\
&= \langle \trans g(e_0 + e_{m+2}), e_0 - e_{m+2} \rangle
\\
&= \langle I_{m+1,2} g^{-1} I_{m+1,2} (e_0 + e_{m+2}), e_0 - e_{m+2} \rangle
\\
&= \langle y, e_0 + e_{m+2} \rangle
\\
&= y_0 + y_{m+2}.
\end{align*}
If $x_0+x_{m+2} \ne 0$,
we set
\begin{subequations}\label{eqn:xy}
  \renewcommand{\theequation}{\theparentequation) (\alph{equation}}%
\abovedisplayskip0pt
\begin{alignat}{3}
&a:=(a_1,\dots,a_{m+1}),
&\quad & a_j := \frac{-y_j}{y_0+y_{m+2}}
&\quad& (1 \le j \le m+1),
\label{eqn:xy(a)}
\\
&b:=(b_1,\dots,b_{m+1}),
&& b_j := \frac{x_j}{x_0+x_{m+2}}
&& (1 \le j \le m+1),
\label{eqn:xy(b)}
\\
&\delta := \operatorname{sgn}(x_0+x_{m+2}),
\hidewidth
\label{eqn:xy(c)}
\\
&t := \log \left| \frac{x_0+x_{m+2}}{2} \right| ,
\hidewidth
\label{eqn:xy(d)}
\\
&m_+ := \delta e^{-tE} \nbar{-b} g\nbar{-a} w_0^{-1}.
\hidewidth
\label{eqn:xy(e)}
\end{alignat}
\end{subequations}

\begin{lemSec}
\label{lem:Bruhat}
Retain the above notation.
\begin{enumerate}
    \renewcommand{\labelenumi}{\upshape\theenumi)}
\item  % 1)
For $g \in O(m+1,2)$, the following three conditions are equivalent:
   \begin{enumerate}
     \renewcommand{\theenumii}{\roman{enumii}}
     \renewcommand{\labelenumii}{\upshape\theenumii)}
   \item  % i)
   $g \notin \overline{\Pmax}$,
   \item  % ii)
   $x_0 + x_{m+2} \ne 0$,
   \item  % iii)
   $y_0 + y_{m+2} \ne 0$.
   \end{enumerate}
\item  % 2)
Suppose one of (therefore, all of) the above conditions holds.
Then, the element $m_+$ defined by \eqref{eqn:xy(e)} belongs to $\Mmax_+
\simeq O(m,1)$.
\end{enumerate}
\end{lemSec}

\subsection{Explicit action of the whole group}
\label{sec:wholeG}

For $g \in O(m+1,2)$,
we set
\begin{align*}
& x \equiv \trans (x_0,\dots,x_{m+2}) := g(e_0 - e_{m+2})
  \quad\text{(see \eqref{eqn:xge}),}
\\
& y \equiv \trans (y_0,\dots,y_{m+2}) := g^{-1}(e_0 - e_{m+2})
  \quad\text{(see \eqref{eqn:yge}).}
\end{align*}
For $g$ such that $x_0 + x_{m+2} \ne 0$,
we also set $a = (a_1,\dots,a_{m+1})$, 
$b = (b_1,\dots,b_{m+1})$, $\delta \in \{ \pm1 \}$,
and $m_+$ as in \eqref{eqn:xy}.

If $x_0 + x_{m+2} = 0$,
then $g \in \overline{\Pmax}$ and $(\pi(g)\psi)(\zeta)$ is obtained
readily by the formulas \eqref{rep_A}--\eqref{rep_N}.
For generic $g$ such that $x_0 + x_{m+2} \ne 0$,
the unitary operator $\pi(g)$ can be written by means of the kernel
function $K(\zeta, \zeta')$ (see \eqref{eqn:KC} for the definition),
the above $x,y \in \mathbb R^{m+3}$ and $m_+ \in \Mmax_+\simeq
O(m,1)$ as follows:

\begin{thmSec}
\label{thm:formulaG}
For $g\in O(m+1,2)$ such that
 $x_0 + x_{m+2} \ne 0$,
the unitary operator $\pi(g)$ is given by the following integral
 formula: for $\psi\in L^2(C)$, 
\begin{align*}
& (\pi(g)\psi)(\zeta)
= \left( \frac{2}{x_0 + x_{m+2}} \right) ^{\frac{m-1}{2}}
  e^{
   \frac{2\sqrt{-1} (x_1 \zeta_1 + \dots + x_{m+1} \zeta_{m+1})}
        {x_0 + x_{m+2} }
         }
\\
&  \int_C K 
   \left( \frac{2\zeta}{|x_0 + x_{m+2} |}, m_+ \zeta' \right)
   e^{
   \frac{-2\sqrt{-1} (y_1 \zeta'_1 + \dots + y_{m+1} \zeta'_{m+1})}
        {y_0 + y_{m+2} }
         }
  \psi (\zeta') du (\zeta').
\end{align*}
\end{thmSec}

\begin{proof}
We write 
$
g = \nbar{b} e^t (\delta m_+) w_0 \nbar{a}
$
as in the form \eqref{eqn:gBruhat}. 
Then, we have
\begin{align*}
&(\pi (g) \psi) (\zeta)
\\
&=\pi(\nbar{b})\pi(e^{tE})\pi(\delta m_+)(\pi(w_0\nbar{a})\psi)
 (\zeta)
\\
&= e^{2\sqrt{-1}\langle b,\zeta\rangle}
  e^{-\frac{m-1}{2} t} \delta^{\frac{m-1}{2}}
     (\pi(w_0\nbar{a})\psi) (e^{-t}\, \trans m_+ \zeta)
  \quad \text{by \eqref{rep_A}--\eqref{rep_N}}.
\\
\intertext{Now, by using Corollary \ref{cor:6.4.1},}
&= e^{2\sqrt{-1}\langle b,\zeta\rangle}
  e^{-\frac{m-1}{2}t} \delta^{\frac{m-1}{2}}
  \int_C K(e^{-t}\, \trans m_+ \zeta, \zeta')
 (\pi(\nbar{a})\psi)(\zeta')du(\zeta') 
\\
&= e^{2\sqrt{-1}\langle b,\zeta\rangle}
  e^{-\frac{m-1}{2} t}
  \delta^{\frac{m-1}{2}}
  \int_C K(e^{-t}\zeta, m_+ \zeta')
  e^{2\sqrt{-1}\langle a,\zeta'\rangle}
  \psi(\zeta')du(\zeta').
\end{align*}
\end{proof}

\section{Appendix: special functions}

For the convenience of the reader, 
we collect here basic formulas of 
special functions in a way that we use 
in this article. 

\subsection{Laguerre polynomials}
\label{Laguerre polynomials}
For $\alpha \in \mathbb C, n \in \mathbb N$, 
the {\it Laguerre polynomials} $L_n^\alpha (x)$ are defined by the formula 
(see \cite[\S 6.2]{xandrews}, for example): 

\begin{alignat}{2}
\label{eqn:defLag}
    L_n^\alpha(x):= {} & \frac{x^{-\alpha}e^{x}}{n!} \frac{d^n}{dx^n}
                      (e^{-x}x^{n+\alpha})   \\
                  = {} & \frac{(\alpha+1)_n}{n!}
                      \sum_{k=0}^n \frac{(-n)_k x^k}{(\alpha+1)_k k!}   
\nonumber
\\
                 = {} & \frac{(-1)^n}{n!} x^n + \cdots + 
                       \frac{(\alpha+1)(\alpha+2)\cdots (\alpha+n)}{n!}.
\nonumber
\end{alignat}
Here, we write $\beta_n$ for $\beta(\beta+1)\cdots(\beta+n-1)$.
The Laguerre polynomial solves the linear ordinary differential equation of second order:
\begin{equation}
\label{ode}
  xu'' + (\alpha+1-x)u'+ nu=0.
\end{equation}
Suppose $\alpha\in\mathbb{R}$ and
 $\alpha >-1$. Then the Laguerre polynomials 
$\set{L_n^\alpha(x)}{n=0,1, \cdots}$ are complete in 
$L^2((0, \infty), x^\alpha e^{-x}dx)$, 
and satisfy the orthogonality relation (see \cite[\S 6.5]{xandrews}):
\begin{equation}
\label{eqn:Lagnorm}
    \int_0^\infty L_m^\alpha(x)L_n^\alpha(x)x^\alpha e^{-x}dx
    = \frac{\Gamma(\alpha+n+1)}{n!}\delta_{mn}  \quad (\alpha>-1).  
\end{equation}

The Hille--Hardy formula gives the bilinear generating function of 
Laguerre polynomials (see \cite[\S10.12 (20)]{xEr}, \cite[\S1 (3)]{xHa}):
\begin{align}
&  \sum_{n=0}^\infty \frac{\Gamma(n+1)}
                                          {\Gamma(n+\alpha+1)}
                                 L_n^\alpha (x) L_n^\alpha (y) w^n   
\nonumber
\\
&  = \frac{1}{1-w} \exp \Bigl( - \frac{(x+y)w}{1-w} \Bigr)
     (-xyw)^{-\frac{\alpha}2} J_\alpha \Bigl( \frac{2\sqrt{-xyw}}{1-w} \Bigr),
  \quad  x>0, y>0, |w|<1.
\label{cor-A}
\end{align}
Here the left-hand side converges absolutely.

\subsection{Hermite polynomials}

Hermite polynomials $H_n(x)$ are given as special values of Laguerre
polynomials. 
We recall from \cite[II, \S10.13]{xEr} that 
\begin{align}
&H_{2m}(x)
= (-1)^m 2^{2m} m!\, L_m^{-\frac{1}{2}} (x^2),
\label{eqn:LHeven}
\\
&H_{2m+1}(x)
= (-1)^m 2^{2m+1} m!\, x L_m^{\frac{1}{2}} (x^2),
\label{eqn:LHodd}
\end{align}
This reduction formula is reflected by the fact that Hermite polynomials
appear in the analysis of the 
Weil representation of
$\SL$,
while Laguerre polynomials appear in the analysis of the minimal
representation of $SO_0(m+1,2)\widetilde{}$
and $\Mp$
(see Remark \ref{rem:4.2.2}).

\subsection{Gegenbauer polynomials}
\label{Gegenbauer polynomials}

For $\nu \in \mathbb C$ and $l \in \mathbb N$, 
the Gegenbauer polynomials $C_l^\nu (x)$ are 
the polynomials of $x$ of degree $l$ given 
 by the {\it Rodrigues formula}
(see \cite[\S 6.4]{xandrews}):
\begin{alignat}{2}
\label{eqn:Gegen}
 C_l^\nu (x) :&=\bigl( -\frac{1}2 \bigr)^l \frac{(2\nu)_l}{l!\, (\nu+\frac{1}2)_l}
                         (1-x^2)^{-\nu+\frac{1}2} \frac{d^l}{dx^l}
                         \bigl( (1-x^2)^{l+\nu-\frac{1}2} \bigr). 
\end{alignat}
It then follows from \eqref{eqn:Gegen} that
\begin{equation*}
C_l^\nu(1)
= \frac{\Gamma(2\nu+l)}{l!\, \Gamma(2\nu)}.
\end{equation*}
We renormalize the Gegenbauer polynomial by
\begin{equation}
\label{eqn:nGegen}
\widetilde{C}_l^\nu(x) := \Gamma(\nu) C_l^\nu(x).
\end{equation}
Then by the duplication formula of the Gamma function:
\begin{equation*}
\Gamma(2\nu) = \frac{2^{2\nu-1}}{\sqrt{\pi}}
               \Gamma(\nu) \Gamma(n+\frac{1}{2}),
\end{equation*}
we have
\begin{equation}
\label{eqn:Gegen1}
\widetilde{C}_l^\nu(1) 
  = \frac{\sqrt{\pi}\, \Gamma(2\nu+l)}
         {2^{2\nu-1} l!\, \Gamma(\nu+\frac{1}{2})}.
\end{equation}
The special value at $\nu=0$ is given by the limit formula
(see \cite[\textbf{I}, \S 3.15.1 (14)]{xEr})
$$
\widetilde{C}_l^0(\cos\theta)
= \lim_{\nu\to0} \Gamma(\nu) C_l^\nu (\cos\theta)
= \frac{2\cos(l\theta)}{l}.
$$
Suppose $\operatorname{Re} \nu > -\frac{1}{2}$.
Then,
$\{\widetilde{C}_l^\nu(x): l=0,1,2,\ldots\}$
forms a complete orthogonal basis of
$L^2((-1,1),(1-x^2)^{\nu-\frac{1}{2}} dx)$,
and the norm of $\widetilde{C}_l^\nu(x)$ is given by
\begin{equation}
\label{eqn:Gegennorm}
\int_{-1}^1 \widetilde{C}_l^\nu(x)^2 (1-x^2)^{\nu-\frac{1}{2}} dx
= \frac{2^{1-2\nu} \pi\Gamma(2\nu+l)}{l!\, (l+\nu)},
\end{equation}
(see \cite[\textbf{I}, \S3.15.1 (17)]{xEr}).
Therefore,
$f\in L^2((-1,1), (1-x^2)^{\nu-\frac{1}{2}} dx)$
has the following expansion:
\begin{equation}
\label{eqn:Gegenexpansion}
f(x) = \sum_{l=0}^\infty \alpha_l^\nu(f) \widetilde{C}_l^\nu(x),
\end{equation}
where we set
\begin{equation}
\label{eqn:coeffGegen}
\alpha_l^\nu(f)
:= \frac{l!\, (l+\nu)}{2^{1-2\nu} \pi\Gamma(2\nu+l)}
\int_{-1}^1 f(x) \widetilde{C}_l^\nu(x) (1-x^2)^{\nu-\frac{1}{2}} dx.
\end{equation}

The following integral formula is used in the proof of Lemma \ref{Lemma5.5}
(see \cite[\S 16.3 (3)]{xerdInt}):
Suppose 
$\operatorname{Re}\beta>-1$ and $\operatorname{Re}\nu>-\frac{1}{2}$.
\begin{multline}\label{eqn:Gegenorth}
 \int_{-1}^1 (1-x)^{\nu-\frac{1}2}(1+x)^\beta \widetilde{C}_n^\nu(x)dx =  \\
  \frac{2^{\beta-\nu+\frac{3}2}\sqrt{\pi}\,\Gamma(\beta+1) 
      \Gamma(2\nu+n)\Gamma(\beta-\nu+\frac{3}2)}
     {n!\, \Gamma(\beta-\nu-n+\frac{3}2)
     \Gamma(\beta+\nu+n+\frac{3}2)}.  
\end{multline}

\subsection{Spherical harmonics and Gegenbauer polynomials}
\label{Spherical harmonics.}
Let $\Delta_{S^{n-1}}$ be the Laplace--Beltrami operator on the $(n-1)$-dimensional
unit sphere $S^{n-1}$. Then the \textit{spherical harmonics} on $S^{n-1}$ are defined as
$$  \mathcal{H}^l(\mathbb{R}^n):=
       \{ f \in C^\infty(S^{n-1}) : \Delta_{S^{n-1}}f=-l(l+n-2)f  \},
\quad
l=0,1,2, \dots. $$
The following facts are well-known:

(1)  $\mathcal{H}^l(\mathbb{R}^n)$ is an irreducible representation
          space of $O(n).$

(2)  It is still irreducible as an $SO(n)$-module if $n\geq 3.$

(3)  $\mathcal{H}^l(\mathbb{R}^2) 
          = \mathbb{C}e^{\sqrt{-1}l\theta} 
            \oplus \mathbb{C}e^{-\sqrt{-1}l\theta},~l\geq 1$ 
           as $SO(2)$-modules, where
           $\theta = \tan^{-1}\frac{x_2}{x_1},~ (x_1, x_2) \in \mathbb{R}^2.$

(4) $\mathcal{H}^l(\mathbb{R}^n)\big|_{O(n-1)}
          \simeq \bigoplus_{k=0}^l \mathcal{H}^k(\mathbb{R}^{n-1})$
          gives an irreducible decomposition of $O(n-1)$-modules. 
          This is also an irreducible decomposition as
          $SO(n-1)$-modules  $(n\geq 4)$.

We regard $O(n-1)$ as the isotropy subgroup of $O(n)$ at 
$(1,0, \cdots, 0) \in \mathbb R^n$. We write $(x_1, \cdots, x_n)$ for the 
standard coordinate of $\mathbb R^n$. Then,
$O(n-1)$-invariant spherical harmonics are unique up to scalar, 
and we have: 
\begin{equation}
\label{eqn:Gegensph}
  \Har{l}{n}^{O(n-1)} \simeq \mathbb C \widetilde{C}_l^\frac{n-2}2 (x_1).
\end{equation}

\subsection{Bessel functions} \label{Bessel functions}

For $\nu \in \mathbb C, z \in \mathbb C \setminus \mathbb R_{\le 0}$, {\it Bessel 
functions} $J_\nu(z)$ are defined by 
\begin{equation}\label{definition of Bessel}
  J_\nu(z):= \bigl(\frac{z}2 \bigr)^\nu
           \sum_{m=0}^\infty \frac{(-1)^m(\frac{z}2)^{2m}}{m!\, \Gamma(\nu+m+1)}.
\end{equation}
It solves the {\it Bessel's differential equation} of second order:
$$
  u''+\frac{1}z u' +\Bigl(1-\frac{\nu^2}{z^2}\Bigr) u=0.
$$
The modified Bessel functions $I_\nu(z)$ are defined by 
\begin{alignat}{2} \label{definition of I-Bessel}
   I_\nu(z):=&\begin{cases} 
                    e^{-\frac{\pi \nu \sqrt{-1}}2} J_\nu(\sqrt{-1}z)  
                  &  (-\pi < \operatorname{arg}z \le \frac{\pi}2),     \\
                    e^{\frac{3 \pi \nu \sqrt{-1}}2} J_\nu(\sqrt{-1}z) 
                  &  (\frac{\pi}2 < \operatorname{arg}z < \pi).
                    \end{cases}   \\
\label{definition of I-Bessel 2}
                 =&\bigl(\frac{z}2 \bigr)^\nu \sum_{n=0}^\infty
                  \frac{\bigl(\frac{z}2\bigr)^{2m}}{m!\, \Gamma(\nu+m+1)}
  \quad (z \in \mathbb C \setminus \mathbb R_{\le 0}).
\end{alignat}
For a special value $\nu=\pm \frac{1}2$, $I_\nu(z)$ reduces to 
\begin{equation}
  I_\frac{1}2 (z)= \sqrt{ \frac{2}{\pi z}} \sinh z, \quad 
  I_{-\frac{1}2}(z)= \sqrt{ \frac{2}{\pi z}} \cosh z.
\end{equation}

We set
\begin{align}
\tilde{J}_\nu(z):={}
&\left(\frac{z}{2}\right)^{-\nu} J_\nu(z)
=
\sum_{m=0}^\infty \frac{(-1)^m \left(\frac{z}{2}\right)^{2m}}
                        {m!\, \Gamma(\nu+m+1)},
\label{eqn:Jtilde}
\\
\tilde{I}_\nu(z):={}
&\left(\frac{z}{2}\right)^{-\nu} I_\nu(z)
=
\sum_{m=0}^\infty \frac{\left(\frac{z}{2}\right)^{2m}}
                        {m!\, \Gamma(\nu+m+1)}.
\label{eqn:Itilde}
\end{align}
We note that $\tilde{J}_\nu(z)$ and $\tilde{I}_\nu(z)$ are entire
functions of $z$, and
$$
\tilde{J}_\nu(\sqrt{-1} z) = \tilde{I}_\nu(z).
$$

The following lemma on the estimate of $I$-Bessel functions are 
used in Subsections \ref{sec:4.3} and \ref{Estimates of the kernel function}.
We need an estimate of $\tilde{I}_\nu(z)$ for $\nu \ge -\frac{1}{2}$:

\begin{lemSec}
\label{4.3.1} 
 There exists a constant $C>0$ such that 
the following estimate holds for all $z \in \mathbb{C}$: 
\begin{alignat}{2}
\label{estimate}
     |\tilde{I}_{\nu} (z)| \leq C  e^{|\operatorname{Re} z|}.
\end{alignat}
\end{lemSec}
\begin{proof}
First, suppose $\nu = -\frac{1}{2}$.
Then
$|\tilde{I}_\nu(z)| = \frac{1}{\sqrt{\pi}}|\cosh z|
 \le \frac{1}{\sqrt\pi} e^{|\operatorname{Re} z|}$.

Next, suppose $\nu > -\frac{1}{2}$.
 By an integral representation of the Bessel function \cite[\S 6.15 (2)]{xWa}:
\begin{equation*}
  \tilde{I}_\nu(z)=\frac{1}{\Gamma(\nu+\frac{1}2)\Gamma(\frac{1}2)}
                  \int_{-1}^1 e^{-zt} (1-t^2)^{\nu-\frac{1}2} dt,
\end{equation*}
we have 
\begin{alignat*}{1}
  |\tilde{I}_\nu(z)| \le  &
                 \frac{1}{\Gamma(\nu+\frac{1}2)\Gamma(\frac{1}2)}
                  \int_{-1}^1 e^{-t \operatorname{Re}z}
                 (1-t^2)^{\operatorname{Re}\nu-\frac{1}2} dt  \\
                 \le &  
                 \frac{e^{|\operatorname{Re}z|}}
                      {\Gamma(\nu+\frac{1}2)\Gamma(\frac{1}2)}
                  \int_{-1}^1 
                  (1-t^2)^{\operatorname{Re}\nu-\frac{1}2} dt
                 \le C e^{|\operatorname{Re} z|}
\end{alignat*}
for some constant $C$ independent of $z$. 
\end{proof} 

The following lemma is used in Subsection \ref{subsec:5.5},
where we set
$\nu= \frac{m-3}2, \alpha = \frac{2 \sqrt{2rr'}}{\sinh \frac{t}2}$.

\begin{lemSec}
\label{Lemma5.5}
Assume $\alpha \in \mathbb C, \nu \ge -\frac{1}2$, and $l \in \mathbb N$. Then 
we have:
\begin{multline}
\label{5.5.1} 
  \int_0^\pi I_\nu\bigl(\alpha \sqrt{1+\cos\theta}\bigr)
  \widetilde{C}_l^{\nu+\frac{1}2} (\cos \theta)
  (1+\cos\theta)^{-\frac{\nu}2}  \sin^{2\nu+1}\theta d\theta  \\
  = \frac{ 2^{\frac{3}2} \sqrt{\pi}\, \Gamma(2\nu+l+1)}
              {\alpha^{\nu+1} l!}
  I_{2\nu+2l+1}\bigl(\sqrt{2}\alpha \bigr).
\end{multline}        
\end{lemSec}

We could not find this formula in the literature, and so we  give its proof here.

\begin{proof}
First we note that the integral  (\ref{5.5.1}) converges
since $I_\nu \bigl(\alpha (\sqrt{1+\cos \theta})\bigr)
\cdot (1+\cos \theta)^{-\frac{\nu}2}$ is continuous on the closed
interval 
$[0, \pi]$. Furthermore, we have a uniformly convergent expansion 
$$
  I_\nu \bigl(\alpha (\sqrt{1+\cos \theta})\bigr)
\cdot (1+\cos \theta)^{-\frac{\nu}2}=
  \sum_{j=0}^\infty  
  \frac{\bigl( \frac{\alpha}2 \bigr)^{\nu+2j} (1+\cos\theta)^j  }
           {j!\, \Gamma(j+\nu+1)}.
$$

Now the left-hand side of (\ref{5.5.1}) equals 
\begin{alignat*}{2}
  &\bigl( \frac{\alpha}2 \bigr)^\nu 
  \sum_{j=0}^\infty 
  \frac{ \bigl( \frac{\alpha}2 \bigr)^{2j} }{j!\, \Gamma(j+\nu+1)}
  \int_0^\pi (1+\cos \theta)^j \widetilde{C}_l^{\nu+\frac{1}2}(\cos \theta)
  \sin^{2\nu+1} \theta d\theta    \\
  &= \bigl( \frac{\alpha}2 \bigr)^\nu 
  \sum_{j=0}^\infty 
  \frac{ \bigl( \frac{\alpha}2 \bigr)^{2j} }{j!\, \Gamma(j+\nu+1)}
  \frac{2^{j+1} \sqrt{\pi}\, \Gamma(j+\nu+1) 
            \Gamma(2\nu+l+1) \Gamma(j+1)}
           {l!\, \Gamma(j-l+1) \Gamma(2\nu+j+l+2)}   \\
  &= \bigl( \frac{\alpha}2 \bigr)^\nu 
  \frac{2 \sqrt{\pi}\, \Gamma(2\nu+l+1)}
           {l!}
  \sum_{j=0}^\infty 
  \frac{ \bigl( \frac{\alpha}{\sqrt 2} \bigr)^{2j}  }
           {\Gamma(j-l+1) \Gamma(2\nu+j+l+2) }    \\
  &= \bigl( \frac{\alpha}2 \bigr)^\nu 
  \frac{2 \sqrt{\pi}\, \Gamma(2\nu+l+1)}
           {l!}
     \bigl( \frac{\alpha}{\sqrt{2}} \bigr)^{2l} 
  \sum_{n=0}^\infty 
  \frac{ \bigl( \frac{\alpha}{\sqrt 2} \bigr)^{n}  }
           {\Gamma(n+1) \Gamma(2\nu+2l+n+2) }    \\
  &= \text{the right-hand side of (\ref{5.5.1}).}
\end{alignat*}
Here the first equality follows from the formula (\ref{eqn:Gegenorth}).
In fact, the substitution of  $x=\cos\theta$ into (\ref{eqn:Gegenorth}) yields
\begin{alignat*}{2}
  &\int_0^\pi (1+\cos \theta)^j \widetilde{C}_l^{\nu+\frac{1}2}(\cos \theta)
  \sin^{2\nu+1} \theta d\theta  \\
  &= \int_{-1}^1 (1-x)^\nu (1+x)^{j+\nu} \widetilde{C}_l^{\nu+\frac{1}2}(x) dx   \\
  &=  \frac{2^{j+1} \sqrt{\pi}\, \Gamma(j+\nu+1) 
            \Gamma(2\nu+l+1) \Gamma(j+1)}
           {l!\, \Gamma(j-l+1) \Gamma(2\nu+j+l+2)}.  
\end{alignat*}

Thus, the lemma has been proved.
\end{proof}

\end{document}